\documentclass[aos,MSNbibl,nameyear,dvips]{arximspdf}
\usepackage{tikz}
\usepackage{graphicx}

\doi{10.1214/14-AOS1206} 
\volume{42}
\issue{4}
\pubyear{2014}
\firstpage{1452}
\lastpage{1482}

\makeatletter

\newcommand{\rrvert}{\vert}
\newcommand{\llvert}{\vert}
\newcommand{\reals}{\mathbb{R}}
\newcommand{\G}{\mathcal{G}}
\newcommand{\X}{\mathfrak{X}}
\newcommand{\vzero}{{0}}
\newcommand{\vone}{{1}}
\newcommand{\implies}{\Longrightarrow}

\newcommand\indep{\protect\mathpalette{\protect\independenT}{\perp}}
\def\independenT#1#2{\mathrel{\rlap{$#1#2$}\mkern2mu{#1#2}}}

\newcommand{\tail}{\operatorname{tail}}
\newcommand{\pa}{\operatorname{pa}}
\newcommand{\an}{\operatorname{an}}
\newcommand{\dec}{\operatorname{de}}
\newcommand{\barren}{\operatorname{barren}}
\newcommand{\dis}{\operatorname{dis}}
\newcommand{\mbl}{\operatorname{mb}}
\newcommand{\pre}{\operatorname{pre}}

\newtheorem{lem}{Lemma}[section]
\newtheorem{teo}[lem]{Theorem}
\newtheorem{prop}[lem]{Proposition}
\newtheorem{cor}[lem]{Corollary}

\newproclaim{dfn}[lem]{Definition}
\newproclaim{rmk}[lem]{Remark}
\newproclaim{exm}[lem]{Example}
\makeatother

\begin{document}
\begin{frontmatter}

\title{Markovian acyclic directed mixed graphs for~discrete data\thanksref{T1}}
\runtitle{Markovian ADMGs for discrete data}

\begin{aug}
\author[a]{\fnms{Robin J.}~\snm{Evans}\corref{}\ead[label=e1]{robin.evans@stats.ox.ac.uk}}
\and
\author[b]{\fnms{Thomas S.} \snm{Richardson}\ead[label=e2]{thomasr@u.washington.edu}}
\runauthor{R. J. Evans and T. S. Richardson}
\affiliation{University of Oxford and University of Washington}
\address[a]{Department of Statistics\\
University of Oxford\\
1--2 South Parks Road\\
Oxford OX1 3TG\\
United Kingdom\\
\printead{e1}}
\address[b]{Department of Statistics\\
University of Washington\\
Box 354322\\
Seattle Washington 98195\\
USA\\
\printead{e2}}
\end{aug}
\thankstext{T1}{Supported by the U.S. National Science Foundation
Grant CNS-0855230 and U.S. National Institutes of Health Grant R01 AI032475.}

\received{\smonth{1} \syear{2013}}
\revised{\smonth{12} \syear{2013}}

%
\begin{abstract}
Acyclic directed mixed graphs (ADMGs) are graphs that
contain directed ($\rightarrow$) and bidirected ($\leftrightarrow$)
edges, subject to the constraint that there are no cycles of directed
edges. Such graphs may be used to represent the conditional
independence structure induced by a DAG model containing hidden
variables on its observed margin. The Markovian model associated with
an ADMG is simply the set of distributions obeying the global Markov
property, given via a simple path criterion (\mbox{m-}separation). We first
present a factorization criterion characterizing the Markovian model
that generalizes the well-known recursive factorization for DAGs. For
the case of finite discrete random variables, we also provide a
parameterization of the model in terms of simple conditional
probabilities, and characterize its variation dependence. We show
that the induced models are smooth. Consequently, Markovian ADMG
models for discrete variables are curved exponential families of
distributions.
\end{abstract}

%
\end{frontmatter}

\section{Introduction}\label{sec1}

A directed graph is a finite collection of vertices, $V$, together
with a collection of ordered pairs $E \subset V \times V$ such that
$(v,v) \notin E$ for any $v$; if $(v,w) \in E$ we write $v \rightarrow
w$. $E$ is the (directed) edge set. We say a directed graph is
\emph{acyclic} if it contains no directed cycles; that is, there is no
sequence of vertices $v_1 \rightarrow v_2 \rightarrow\cdots
\rightarrow v_k \rightarrow v_1$, for any $k > 1$. We call such a
graph a \emph{directed acyclic graph} (DAG). Models based on DAGs are
popular because
of their simple definition in terms of a recursive factorization, easy
to determine conditional independence constraints, and potential for
causal interpretations [\citeauthor{pearlbiom} (\citeyear{pearlbiom,pearl2009}) \citet{cps93,robinsmcm2011}]. Unfortunately, if some of
the variables in a DAG
are unobserved, the resulting pattern of conditional independences no
longer corresponds to a DAG model (on the observed variables); in this
sense, DAGs are not closed under marginalization.

\begin{figure}

\includegraphics{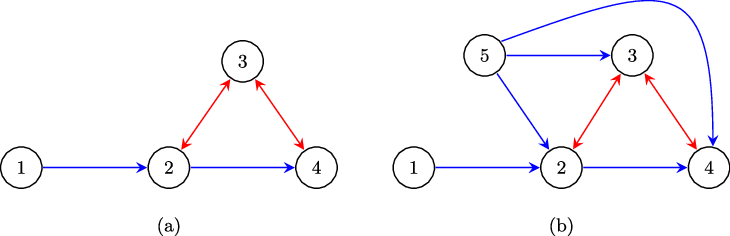}

\caption{\textup{(a)}~An acyclic directed mixed graph, $\mathcal{L}$.
\textup{(b)}~An ADMG studied by \citet{evansrichardson2013}.}\label{figgraph1}
\end{figure}

An \emph{acyclic directed mixed graph} (ADMG) consists of a DAG with
vertices $V$ and directed edges $E$, together with a collection $B$ of
unordered (distinct) pairs of elements of $V$; these are the \emph{bidirected edges}. If $\{v,w\} \in B$ we
write $v \leftrightarrow w$, and if in addition $(v,w) \in E$ this is
denoted
$v\; {\tikz[baseline=-2pt]{
\draw[<->] (0cm,0cm) -- (3mm,0cm);
\draw[->] (0cm,4pt) -- (3mm,4pt);}
}\; w$.
Graphical definitions are best understood visually, so we invite the
reader to consult the example ADMGs given in Figure~\ref{figgraph1}.

%

Like DAGs, acyclic directed mixed graphs can be interpreted, via a
Markov property, as representing a set of probability distributions
defined by conditional independence restrictions; these can be read
off the graph using a graphical separation criterion. The advantage
of ADMGs is that they are closed under marginalization, in the
sense mentioned above [\citet{richardson02}]; indeed they represent
precisely the conditional independence relations which can be obtained
by marginalizing DAGs. \citet{richardson03} gave a global Markov
property and ordered local Markov property for ADMG models, and showed
their equivalence.

The patterns of conditional independence implied by a DAG give rise to
curved exponential families in the case of discrete random variables
and, therefore, these models have well understood asymptotic statistical
properties. However, general models induced by conditional
independence constraints do not share this property, and it may be
challenging to determine their dimension; for example, certain
interpretations of chain graphs are known to lead to non-smooth models
[\citet{drton09}]. In this paper, we show that discrete ADMG
models are
curved exponential families, and give a smooth parameterization.

\citet{evansrichardson2013} provide a number of applied examples for
\mbox{ADMGs} representing discrete distributions---such as using the graph in
Figure~\ref{figgraph1}(b) to model an encouragement design for an influenza
vaccine---and they discuss the relationship between Markovian \mbox{ADMG} models
and marginal log-linear models [\citet{br02,bart07}].
\mbox{ADMGs} also arise in studying
general conditions for identifying intervention distributions, under
the causal interpretation of a DAG model
[see \citet{dawiddidelez2010,huang06do,pearlrobins1995}, \citeauthor{shpitser06id} (\citeyear{shpitser06id,shpitser06idc}), \citet{tian02b,silva09}].

%
\begin{figure}

\includegraphics{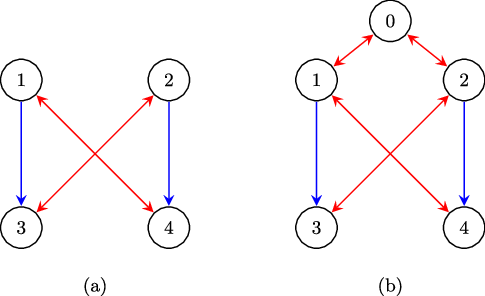}

\caption{\textup{(a)}~An ADMG in which there is no vertex ordering such that
all parents of a head precede every vertex in the head;
\textup{(b)} $\{0,3,4\}$ forms a head in this ADMG, but the induced subgraph on
$\{0,3,4\}$ is not connected.}\label{fignoorder}
\end{figure}

This paper provides a factorization criterion for joint distributions
obeying the global Markov property with respect to an ADMG as well as
a parameterization of these models in the discrete case.
The factorizations so obtained are unusual: the graph
in Figure~\ref{fignoorder}(a), for example, gives
\begin{eqnarray*}
f_{1234}(x_1, x_2, x_3,
x_4) &=& f_{23|1}(x_2, x_3 \mid
x_1) \cdot f_{14|2}(x_1, x_4 \mid
x_2),
\end{eqnarray*}
showing that the joint distribution is a product of two conditional
distributions that we would not usually expect to multiply together
(see Example~\ref{exmfactor1}). The factorization criterion
generalizes the well known one for DAGs, and is analogous to the
Hammersley--Clifford theorem for undirected graphical models
[\citet{hammersley71}]; the parameterization enables model
fitting, and
is used to prove that the discrete models are curved exponential
families of distributions.

ADMGs may be viewed as a subclass of the larger classes of
summary graphs [\citet{wermuth11}] and ribbonless mixed graphs
[\citet{sadeghi11b,sadeghi11}], which allow for undirected edges. The
factorization and parameterization developed here may be extended to
these larger classes without difficulty.

The remainder of the paper is organized as follows: Section~\ref{secgraphical} introduces basic graphical concepts. In Section~\ref
{secpartitions}, we give conditions under which a partial ordering
on a class of subsets may be used to define partitions of arbitrary
subsets. In Section~\ref{secfact}, we use these tools to develop our
factorization criterion, which then forms the basis of the simple
parameterization introduced in Section~\ref{sectpara}. In Section~\ref
{secsmooth}, we show that the Markov model associated with an ADMG
is smooth, and characterize the variation dependence of the
parameterization. Finally, Section~\ref{secdiscussion} contains a
brief discussion.

\section{Graphical definitions and Markov properties}\label{secgraphical}

Let $\G$ be an acyclic directed mixed graph with vertices $V$; the
\emph{induced subgraph} of $\G$ over $A \subseteq V$, denoted $\G_A$,
is the graph with vertex set $A$, and all those (directed or bidirected) edges which join two
vertices that are both in $A$.

A \emph{path} in $\G$ is a sequence of adjacent edges, without
repetition of a vertex; a~path may be empty, or equivalently consist
of only one vertex. The first and last vertices on a path are the
\emph{endpoints} (these are not distinct if the path is empty); other
vertices on the path (if any) are \emph{non-endpoints}. The graph
$\mathcal{L}$ in Figure~\ref{figgraph1}(a), for example, contains the
path $1 \rightarrow2 \rightarrow4 \leftrightarrow3$, with endpoints
1 and 3, and non-endpoints 2 and 4. A \emph{directed} path is one in
which all the edges are directed ($\rightarrow$) and are oriented in
the same direction, whereas a \emph{bidirected path} consists entirely
of bidirected edges.

We use the usual familial terminology for vertices in a graph. If $w
\rightarrow v$, we say that $w$ is a \emph{parent} of $v$; the set of
parents of $v$ is denoted $\pa_\G(v)$. More generally, $w$ is an
\emph{ancestor} of $v$ if there is a directed path from $w$ to $v$
(note that this includes the case $v=w$); conversely $v$ is a
\emph{descendant} of $w$. The ancestors and descendants of $v$ are
denoted $\an_\G(v)$ and $\dec_\G(v)$, respectively. In the graph
$\mathcal{L}$ in Figure~\ref{figgraph1}(a), for instance, the
ancestors of 4 are the vertices $\an_{\mathcal{L}}(4) = \{1,2,4\}$,
and
\[
\pa_{\mathcal{L}}(4) = \{2\}, \qquad\dec_{\mathcal{L}}(4) = \{4\}.
\]
The \emph{district} containing $v$, denoted $\dis_\G(v)$, is the set
of vertices $w$ such that $v \leftrightarrow\cdots\leftrightarrow
w$, including $v$ itself; for example, the district of 4 in
$\mathcal{L}$ is $\{2,3,4\}$.
We apply these functions disjunctively to sets so that, for example,
\[
\an_\G(W) = \bigcup_{v \in W}
\an_\G(v).
\]

A set of vertices $A$ is \emph{ancestral} if $A = \an_{\G} (A)$; that
is, $A$ contains all its own ancestors. Define
\[
\barren_\G(B) \equiv \bigl\{v \in B | \dec_\G(v) \cap B =
\{v\}\bigr\}.
\]
We say a set $B$ is \emph{barren} if $B = \barren_\G(B)$; that is, it
contains none of its nontrivial descendants in $\G$. We will also
use the notation $\dis_A(v)$ as a shorthand for $\dis_{\G_A}(v)$, the
{\em district containing $v$ in the induced subgraph of $\G$ on
$A$}.\vadjust{\goodbreak}

For an ADMG $\G$ with vertex set $V$, we consider collections of
random variables $(X_{v})_{ v\in V}$ taking values in probability
spaces $({\X}_{v})_{v\in V}$; these spaces are either finite discrete
sets or finite-dimensional real vector spaces. For $A\subseteq V$, we
let ${\mathfrak X}_{A}\equiv
\mbox{\fontsize{18}{18}\selectfont{$\times$}}_{v\in A}({\X
}_{v})$, ${\X}\equiv
{\X}_{V}$ and $X_{A}\equiv(X_{v})_{v\in A}$. We abuse notation in the
usual way: $v$~denotes both a vertex and the random variable $X_{v}$,
likewise $A$ denotes both a set of vertices and the random vector
$X_{A}$. For fixed elements of $\X_v$ and $\X_A$, we write $x_v$ and
$x_A$, respectively.

The relationship between a graph $\G$ and random variables $X_V$ is
governed by Markov properties specified in terms of paths.
A non-endpoint vertex $c$ on a path~$\pi$, is a \emph{collider} on $\pi$
if the edges preceding and succeeding $c$ on the path both
 have an arrowhead
at $c$, for example, $\rightarrow c \leftarrow$ or $\leftrightarrow
c\leftarrow$; otherwise $c$ is a \emph{non-collider}.

%
\begin{dfn}
A path $\pi$ in $\G$ between two vertices $v, w \in V(\G)$ is said to
be \emph{blocked} by a set $C \subseteq V\setminus\{v,w\}$ if either:
\begin{longlist}[(ii)]
\item there is a non-collider on $\pi$, and that non-collider is
contained in $C$; or
\item there is a collider on $\pi$ which is not in $\an_\G(C)$.
\end{longlist}
We say $v$ and $w$ are \emph{\mbox{m-}separated} given $C$ in $\G$ if every
path from $v$ to $w$ in $\G$ is blocked by $C$. Note that $C$ may be empty.
Sets $A, B \subseteq V$ are said to be \mbox{m-}separated given $C \subseteq
V\setminus(A \cup B)$ if every pair $a \in A$ and $b \in B$ are
\mbox{m-}separated given $C$.
\end{dfn}

The special case of \mbox{m-}separation for DAGs is the better known
d-separation [\citet{pearl88,lau96}]. We next relate \mbox{m-}separation
to conditional independence, for which we use the now standard notation
of \citet{dawidcondind}: for random variables $X$, $Y$ and $Z$ we
denote the statement ``$X$ is independent of $Y$ conditional on $Z$'' by
$X \indep Y | Z$. If $Z$ is empty, we write $X \indep Y$.

%
\begin{dfn}
A probability measure $P$ on $\mathfrak X$ is said to satisfy the \emph
{global Markov property} (GMP) for an acyclic directed mixed graph $\G
$, if for all disjoint sets $A, B, C \subseteq V$ with $A$ and $B$
nonempty, $A$ being \mbox{m-}separated from $B$ given $C$ implies that $X_A
\indep X_B | X_C$ under $P$.
\end{dfn}

Consider the ADMG $\mathcal{L}$ in Figure~\ref{figgraph1}(a); the
vertices 1 and 4 are \mbox{m-}separated conditional on 2, and 1 and 3 are
\mbox{m-}separated unconditionally. It is not hard to verify that no other
\mbox{m-}separation relations hold for this graph, and that therefore a
distribution $P$ obeys the global Markov property with respect to $\G$
if and only if $X_1 \indep X_4 | X_2$ and $X_1 \indep X_3$ under $P$.

%
\begin{dfn} \label{dfnblanket}
Let $\G$ be an ADMG containing an ancestral set $A$ and a vertex
$v \in\barren_\G(A)$. Define
\[
\mbl_\G(v, A) \equiv\pa_\G\bigl( \dis_A
(v) \bigr) \cup\bigl( \dis_A (v) \setminus\{v\} \bigr)
\]
to be the \emph{Markov blanket} for $v$ in the induced subgraph
$\G_A$.
For a set of vertices $W \subseteq\barren_\G(A)$, we analogously
define the Markov blanket of $W$ to be
\[
\mbl_\G(W, A) \equiv\pa_\G\bigl( \dis_A
(W) \bigr) \cup\bigl( \dis_A (W) \setminus W \bigr).
\]
Let $<$ be a \emph{topological ordering} on the vertices of $\G$,
meaning that no vertex appears before any of its ancestors; let
$\pre_{\G,<} (v)$ be the set of vertices containing $v$ and all
vertices preceding $v$ in the ordering. A probability measure $P$ is
said to satisfy the \emph{ordered local Markov property} for $\G$ with
respect to $<$, if for any $v$ and ancestral set $A$ such that $v \in
A \subseteq\pre_{\G,<} (v)$,
\[
v \indep A \setminus\bigl(\mbl_{\G}(v, A) \cup\{v\}\bigr) |
\mbl_{\G}(v, A)
\]
with respect to $P$.
\end{dfn}

%
\begin{rmk}
For $v \in\barren_\G(A)$, the Markov blanket for $v$ in $A$
consists of those vertices in $A \setminus\{v\}$ that can be
reached from $v$ by paths through $A$ on which all non-endpoints are colliders.
\end{rmk}

%
\begin{exm}
One can easily verify that for the graph in Figure~\ref{figgraph1}(a),
\[
\mbl_{\mathcal{L}}\bigl(4, \{1,2,4\}\bigr) = \{2\}, \qquad\mbl
_{\mathcal{L}}
\bigl(3, \{1,3\}\bigr) = \varnothing,
\]
and that therefore under the topological ordering $1,2,3,4$, the
ordered local Markov property implies $X_4 \indep X_1 | X_2$ and
$X_3 \indep X_1$, just as the global Markov property does.
\end{exm}

The following result shows that the two properties are, in fact,
always equivalent.

%
\begin{prop}[{[\citet{richardson03}, Theorem~2]}] \label{propordered}
Let $\G$ be an ADMG, and $<$ a topological ordering of its vertices;
further let $P$ be a probability measure on~$\X_V$. The following are
equivalent:
\begin{longlist}[(ii)]
\item[(i)] $P$ obeys the global Markov property with respect to $\G$;
\item[(ii)] $P$ obeys the ordered local Markov property with respect to $\G$
and $<$.
\end{longlist}
\end{prop}

In particular, this result implies that if the ordered local Markov
property is satisfied for some topological ordering $<$, then it is
satisfied for all such orderings.


\section{Partitions and partial orderings} \label{secpartitions}

The global Markov property for DAGs can be equivalently stated in terms
of a simple factorization criterion applied to the joint distribution.
In order to achieve something similar for ADMGs, we will need to
consider partitions of sets of vertices into appropriate blocks. This
section develops the necessary mathematical theory on functions that
define partitions.\vadjust{\goodbreak}

Let $V$ be an arbitrary finite set, and let $\mathcal{H}$ be a
collection of nonempty subsets of $V$, with the restriction that $\{v\}
\in\mathcal{H}$ for all $v \in V$ (i.e., all singletons are in
$\mathcal{H}$). Let $\prec$ be a partial ordering on the elements of
$\mathcal{H}$, and write $H_1 \preceq H_2$ to mean that either $H_1
\prec H_2$ or $H_1 = H_2$.

%
\begin{dfn}
We say that $\prec$ is \emph{partition-suitable} (for $\mathcal{H}$)
if for any $H_1, H_2 \in\mathcal{H}$ with $H_1 \cap H_2 \neq
\varnothing$, there exists $H^* \in\mathcal{H}$ such that $H^*
\subseteq H_1 \cup H_2$ and $H_i \preceq H^*$ for each $i=1,2$.
\end{dfn}

In other words, partition-suitability requires that any two
intersecting elements of $\mathcal{H}$ are dominated with respect to
$\prec$ by some element of $\mathcal{H}$.

Define a function $\Phi$ on subsets of $V$ such that $\Phi(W)$ ``picks
out'' the $\prec$-maximal elements of $\mathcal{H}$ which are subsets
of $W$. That is, it returns the collection of subsets
\[
\Phi(W) \equiv\bigl\{H \in\mathcal{H} | H \subseteq W\mbox{ and } H
\nprec H'\mbox{ for all other } H' \subseteq W\bigr\}.
\]
Partition-suitability ensures that the sets in $\Phi(W)$ are disjoint.

%
\begin{prop} \label{propdisjoint}
If $\prec$ is partition-suitable and $H_1, H_2 \in \Phi(A)$ for some
set~$A$, then either $H_1 = H_2$ or $H_1 \cap H_2 = \varnothing$.
\end{prop}

\begin{pf}
This is immediate from the definition of partition-suitable.
\end{pf}

Now let
\[
\psi(W) \equiv W \setminus\bigcup_{C \in\Phi(W)} C,
\]
that is, $\psi$ returns those elements of $W$ which are not contained
in any set in $\Phi(W)$.
Then recursively define a \emph{partitioning function} $[\cdot]$ on
subsets of $V$ by $[\varnothing] = \varnothing$, and
\[
[W] \equiv\Phi(W) \cup\bigl[\psi(W)\bigr].
\]
The idea is that the function $\Phi$ ``removes'' the maximal sets from
$W$, and the procedure is then applied again to what remains,
$\psi(W)$. The following proposition shows that each vertex of $W$ is
contained within precisely one set in $[W]$.

%
\begin{prop}
If $\prec$ is partition-suitable, then the function $[\cdot]$
partitions sets. That is, for any $W \subseteq V$,
\[
\bigcup_{H \in[W]} H = W,
\]
and if $A, B \in[W]$ then either $A=B$ or $A \cap B = \varnothing$.
\end{prop}

\begin{pf}
We proceed by induction on the size of $W$. If $W = \varnothing$ the
result follows from the definition. Also by definition, if $W \neq
\varnothing$ then
\[
[W] = \Phi(W) \cup\bigl[\psi(W)\bigr],
\]
so the induction hypothesis and the definitions of $\Phi$ and $\psi$
mean we need only check that $\Phi(W)$ is nonempty and contains
disjoint sets.

The first claim follows from the fact that $\prec$ is a partial
ordering, and so always contains at least one maximal element (since
$V$ is finite); the second is a direct application of Proposition
\ref{propdisjoint}.
\end{pf}

%
\begin{lem} \label{lemheadfirst}
Let $\prec$ be partition-suitable, $A \subseteq V$ and $H \in
\Phi(A)$. If $H \subseteq B \subseteq A$ for some subset $B$, then
$H \in\Phi(B)$.
\end{lem}

\begin{pf}
Let $\mathcal{H}_A$ be the set of subsets in $\mathcal{H}$ contained
within $A$. If $H \in\Phi(A) \subseteq\mathcal{H}_A$ then $H$ is
maximal with respect to $\prec$ in $\mathcal{H}_A$. It is trivial
that $\mathcal{H}_B \subseteq\mathcal{H}_A$, and so $H$ is also
maximal in $\mathcal{H}_B$. Thus, $H \in\Phi(B)$.
\end{pf}

We can paraphrase Lemma \ref{lemheadfirst} as saying that if a set
$H$ is removed from $A$ at the first application of $\Phi$, then $H$
is contained in the partition of any subset $B$ of~$A$ (provided $B$
contains $H$).

The next proposition shows that partitioning functions as we have
defined them are stable when some set in the partition is removed.
This ``stability'' is very useful when trying to understand
the properties of the partition.

%
\begin{prop} \label{proprmv}
If $C \in[W]$, then $[W] = \{C\} \cup[W \setminus C]$.
\end{prop}

\begin{pf}
We proceed by induction on the size of $W$. If $[W] = \{C\}$,
including any case in which $|W|=1$, the result is trivial.

If $C$ is not maximal with respect to $\prec$ in $W$ then, by Lemma
\ref{lemheadfirst}, $\Phi(W) = \Phi(W \setminus C)$, so
\begin{eqnarray*}
[W] &=& \Phi(W) \cup\bigl[\psi(W)\bigr]
\\
&=& \Phi(W \setminus C) \cup\bigl[\psi(W)\bigr],
\end{eqnarray*}
and the problem reduces to showing that
\[
\bigl[\psi(W)\bigr] = \{C\} \cup\bigl[\psi(W \setminus C)\bigr] = \{
C\} \cup
\bigl[\psi(W) \setminus C\bigr],
\]
which holds by the induction hypothesis. Thus, without loss of
generality, suppose $C \in\Phi(W)$.

Now, by Lemma \ref{lemheadfirst} and the supposition, $\Phi(W
\setminus C) \cup\{C\} \supseteq\Phi(W)$, and if equality holds we
are done. Otherwise let $C_1, \ldots, C_k$ be the sets in $\Phi(W
\setminus C)$ but not in $\Phi(W)$. Note that by definition, $C_1,
\ldots, C_k \subseteq\psi(W)$. Further, these sets are maximal in $W
\setminus C$, so by Lemma \ref{lemheadfirst} they are also maximal in
$\psi(W) \subseteq W \setminus C$. Then the problem reduces to
showing that
\[
\bigl[\psi(W)\bigr] = \{C_1, \ldots, C_k\} \cup\bigl[
\psi(W) \setminus(C_1 \cup\cdots\cup C_k)\bigr],
\]
which follows from repeated application of the induction hypothesis.
\end{pf}

Lastly, we show that if each set in $\mathcal{H}$ is contained within a piece of
some partition of $V$, then the partitioning function can
be applied separately to each piece of this coarser partition.

%
\begin{prop} \label{proppart}
Let $D_1, \ldots, D_k$ be a partition of $V$, and suppose that every $H
\in\mathcal{H}$ is contained within some $D_i$. Let $\prec$ be a
partition-suitable partial ordering on $\mathcal{H}$. Then for all $W
\subseteq V$,
\[
[W] = \bigcup_{i=1}^k [W \cap
D_i].
\]
\end{prop}

\begin{pf}
We prove the case $k=2$, from which the general result follows by
repeated applications. If either of $W \cap D_1$ or $W \cap D_2$ are
empty, then the result is trivial. By definitions
\[
[W] = \Phi(W) \cup\bigl[\psi(W)\bigr];
\]
$\psi(W)$ is strictly smaller than $W$, so by the induction hypothesis
\[
[W] = \Phi(W) \cup\bigl[\psi(W) \cap D_1\bigr] \cup\bigl[\psi(W)
\cap D_2\bigr].
\]
Define $\mathcal{C}_1$, $\mathcal{C}_2$ so that $\Phi(W) =
\mathcal{C}_1 \cup\mathcal{C}_2$ and each $H \in\mathcal{C}_i$ is a
subset of $D_i$ only; since the elements of $\mathcal{C}_i$ are
maximal with respect to $\prec$ in $W$, by Lemma \ref{lemheadfirst}
they are also maximal in $W \cap D_i$. Hence, $\mathcal{C}_i \subseteq
\Phi(W \cap D_i)$. Repeatedly applying Proposition
\ref{proprmv} gives
%
\begin{eqnarray*}
\mathcal{C}_i \cup\bigl[\psi(W) \cap D_i\bigr]
&=& [W \cap D_i],
\end{eqnarray*}
because $(\psi(W) \cap D_i) \cup\bigcup_{C \in\mathcal{C}_i} C = W
\cap D_i$.
Hence the result.
\end{pf}


\section{The factorization criterion} \label{secfact}

Let $P$ be a probability measure having density $f_V\dvtx \X_V
\rightarrow\reals$ with respect to some $\sigma$-finite dominating
product measure $\mu$ on~$\X_V$. For $U,W \subseteq V$, we denote by
$f_W\dvtx \X_W \rightarrow\reals$ the marginal density over $W$, and by
$f_{W | U}(\cdot | u)\dvtx \X_W \rightarrow\reals$ for $f_U(u) > 0$
the conditional density of $W$ given $U=u$ (more precisely: any member
of the equivalence class of such densities). Then $P$ obeys the
global Markov property with respect to a DAG if and only if it
factorizes as
\[
f_V(x_V) = \prod_{v \in V}
f_{v|\pa(v)}(x_v | x_{\pa(v)}),
\]
for $\mu$-almost all $x_V \in\X_V$ [see, e.g., \citet{lau96}].
In the sequel, all equalities over $f$ are considered to
hold almost everywhere with respect to~$\mu$.

In this section, we show that factorizations can also be used to
characterize Markov models over ADMGs; however, as we shall see, the
criterion is more complicated than that for DAGs.

%
\begin{exm}
Consider the ADMG in Figure~\ref{figgraph1}(a).
A distribution which obeys the global Markov property with respect to
this graph satisfies $X_1 \indep X_3$ and $X_1 \indep X_4 | X_2$. It is
not possible to specify a factorization on the joint distribution of
$X_1$, $X_2$, $X_3$ and $X_4$ which implies precisely these two
independences. Instead, we require factorizations of certain marginal
distributions:
\begin{eqnarray*}
f_{13}(x_1, x_3) &=& f_1(x_1)
\cdot f_3(x_3),
\\
f_{124}(x_1, x_2, x_4) &=&
f_1(x_1) \cdot f_{2|1}(x_2 |
x_1) \cdot f_{4|2}(x_4 | x_2).
\end{eqnarray*}
\end{exm}

Such marginal factorizations can be
used to represent distributions which obey the global Markov property
with respect to an ADMG.

%
\begin{dfn}[(Head)]
%
%
A vertex set $H \subseteq V$ is a \emph{head} if it is barren in $\G$
and  contained within a single district of $\G_{\an(H)}$. We
write $\mathcal{H}(\G)$ for the collection of all heads in $\G$.
\end{dfn}

Note that every singleton vertex $\{v\}$ forms a head.

%
\begin{exm}
For the ADMG shown in Figure~\ref{fignoorder}(b), we have the following:
\begin{eqnarray*}
\mathcal{H}(\G) &=& \bigl\{ \{0\}, \{1\}, \{2\}, \{3\}, \{4\}, \{0,1\},
\{0,2\},
\{1,4\}, \{2,3\},
\\
&&\hspace*{59pt} \{0,1,2\}, \{0,1,4\}, \{0,2,3\}, \{0,3,4\} \bigr\}.
\end{eqnarray*}
Notice that although they are contained within a single district, the
sets $\{0,1,2,4\}$, $\{0,1,2,3\}$ and $\{0,1,2,3,4\}$ do not form heads
because they are not barren.
Also observe that $\{0,3,4\}$ does form a head, even though the induced
subgraph $\G_{\{0,3,4\}}$ is not connected [because
$\{0,3,4\}$ is a subset of a single district in $\G_{\an(\{0,3,4\})}$,
as required].
\end{exm}

%
\begin{dfn}[(Tail)]
For any head $H$, the \emph{tail} of $H$ is the set
\[
\tail_\G(H) \equiv\bigl(\dis_{\an(H)}(H) \setminus H\bigr)
\cup\pa\bigl(\dis_{\an(H)}(H)\bigr).
\]
%
If the context makes it clear which head we are referring to, we will
sometimes denote a tail simply by $T$.
\end{dfn}

Note that the tail is a subset of the ancestors of the head. An
intuitive interpretation is that a head $H$ is a
set within which no independence relations hold without marginalizing
some elements of $H$, and the tail is the Markov blanket for $H$
within the set $\an_\G(H)$. We can therefore factorize ancestral sets
into heads conditional upon their tail sets; see Remark \ref{rmkmbl}
below.

%
\begin{exm}
In the special case of a DAG, the heads are precisely all singleton
vertices $\{v\}$, and the tails are the sets of parents $\pa_{\G}(v)$.
In a purely bidirected graph, the heads are just the connected sets,
and the tails are all empty.
\end{exm}

%
\begin{exm}
The graph $\mathcal{L}$ in Figure~\ref{figgraph1}(a) has the following
head--tail pairs:\vspace*{6pt}
\begin{center}
\begin{tabular}{@{}lcccccc@{}}
\hline
$H$&$\{1\}$&$\{2\}$&$\{3\}$&$\{2,3\}$&$\{4\}$&$\{3,4\}$\\
$T$&$\varnothing$&$\{1\}$&$\varnothing$&$\{1\}$&$\{2\}$&$\{1,2\}$\\
\hline
\end{tabular}\vspace*{6pt}
\end{center}
Note that the set $\{2,3,4\}$ is not a head, because it is not barren.
\end{exm}



In general, it is not possible to order the vertices in an acyclic
directed mixed graph such that, for each head $H$, all the vertices
in $\pa_\G(H)$ precede all the vertices in $H$. A counterexample is
given in Figure~\ref{fignoorder}(a), which is taken from
\citet{richardson09}. The head $\{1,4\}$ has parent 2, and the head
$\{2,3\}$ has parent 1, so whichever way we order the vertices 1 and
2, the condition will be violated.

However, there is a well-defined partial ordering on heads which will
be useful to us, and satisfies the essential property of
partition-suitability from Section~\ref{secpartitions}.

%
\begin{dfn} \label{dfnprec}
For two distinct heads $H_i$ and $H_j$ in an ADMG $\G$, say that $H_i
\prec H_j$ if $H_i \subseteq\an_\G(H_j)$.
\end{dfn}

%
\begin{lem} \label{lemprec}
The (strict) partial ordering $\prec$ is well defined.
\end{lem}

\begin{pf}
We need to verify that $\prec$ is irreflexive, asymmetric and
transitive; irreflexivity is by definition. Asymmetry amounts to
$H_i \prec H_j \implies H_j \nprec H_i$; suppose not for contradiction,
so that there exist distinct heads $H_i$ and $H_j$ with $H_i \prec H_j$
and $H_j \prec H_i$. Since $H_i$ and $H_j$ are distinct, there exists a
vertex $v$ which is in one of these heads but not the other; assume
without loss of generality that $v \in H_j \setminus H_i$.

Since $H_j \subseteq\an_\G(H_i)$, we can find a directed path $\pi_1$
from $v$ to some vertex $w \in H_i$; the path is nonempty because $v
\notin H_i$. However, since we also have $H_i \subseteq\an_\G(H_j)$,
we can find a (possibly empty) directed path $\pi_2$ from $w$ to some
$x \in H_j$. Now, the concatenation of $\pi_1$ and $\pi_2$ is also a
path, because any repeated vertices would imply a directed cycle in the
graph. Call this new path $\pi$.

But $\pi$ is a nonempty directed path between two vertices in $H_j$,
which violates the requirement that heads are barren. Hence, asymmetry holds.

For transitivity, if $H_i \prec H_j$ and $H_j \prec H_k$, then clearly
we can find a directed path from any element $v \in H_i$ to some
element of $H_k$, simply by concatenating paths from $v \in H_i$ to
some $w \in H_j$ and from $w$ to $H_k$. Hence, $H_i \subseteq\an_\G
(H_k)$, and so $H_i \prec H_k$.
\end{pf}


%
\begin{lem} \label{lempartsuit}
The partial ordering $\prec$ on the heads $\mathcal{H}(\G)$ of an ADMG
$\G$ is partition-suitable.
\end{lem}

\begin{pf*}{Proof sketch (see the \hyperref[app]{Appendix} for details)}
If two heads $H_1, H_2$ are distinct and $H_1 \cap H_2 \neq\varnothing
$, then
$H^* = \barren_\G(H_1 \cup H_2)$ is a head, $H_1 \preceq H^*$ and $H_2
\preceq H^*$.
\end{pf*}

Note that in general $H^*$ may be a strict subset of $H_1 \cup
H_2$. For example, consider the graph shown in Figure~\ref
{fignoorder}(b), and let $H_1 = \{0,1,4\}$ and $H_2 =
\{0,2,3\}$ so that $H_1, H_2 \in{\mathcal H}(\G)$ and $H_1\cap H_2
=\{0\}$. However, $H^* = \barren_\G(H_1 \cup H_2) = \{0,3,4\}
\subsetneq H_1\cup H_2$.

Denote the relevant functions from Section~\ref{secpartitions}
defined by this partial ordering by $\Phi_\G$, $\psi_\G$ and
$[\cdot]_\G$, respectively. This partitioning function allows us to
factorize probabilities for ADMGs into expressions based upon heads
and tails.

\begin{exm}
For the graph $\mathcal{L}$ in Figure~\ref{figgraph1}(a), we have\vspace*{6pt}
\begin{center}
\begin{tabular}{@{}lcccccc@{}}
\hline
$H$&$\{1\}$&$\{2\}$&$\{3\}$&$\{2,3\}$&$\{4\}$&$\{3, 4\}$\\
$\an_\G(H)$&$\{1\}$&$\{1,2\}$&$\{3\}$&$\{1,2,3\}$&$\{1,2,4\}$&$\{
1,2,3,4\}$\\
\hline
\end{tabular}
\end{center}\vspace*{6pt}
so that
\begin{eqnarray*}
\{1\} &\prec&\{2\} \prec\{2,3\} \prec\{3,4\},
\\
\{2\} &\prec& \{4\} \prec\{3,4\}, \qquad\{3\} \prec\{2,3\}.
\end{eqnarray*}
Then, for example, $\Phi_{\mathcal{L}}(\{2,3,4\}) = \{\{3,4\}\}$, and
$\Phi_{\mathcal{L}}(\{2\}) = \{\{2\}\}$, giving
\[
\bigl[\{2,3,4\}\bigr]_{\mathcal{L}} = \bigl\{\{3,4\},\{2\}\bigr\}.
\]
\end{exm}

%
\begin{exm}
For the graph in Figure~\ref{fignoorder}(a), we have\vspace*{6pt}
\begin{center}
\begin{tabular}{@{}lcccccc@{}}
\hline
$H$&$\{1\}$&$\{2\}$ & $\{3\}$& $\{4\}$&$\{1,4\}$&$\{2,3\}$\\
$\an_\G(H)$&$\{1\}$&$\{2\}$&$\{1,3\}$&$\{2,4\}$&$\{1,2,4\}$&$\{1,2,3\}$\\
\hline
\end{tabular}\vspace*{6pt}
\end{center}
Thus
$\{1\} \prec\{3\} \prec\{2,3\}\succ\{2\}$ and $\{2\} \prec\{4\}\prec\{1,4\}\succ \{1\}$.
\end{exm}

Now we can provide a factorization criterion for acyclic directed mixed graphs.

%
\begin{teo} \label{teofactorization}
Let $\G$ be an ADMG, and $P$ a probability distribution on $\X_V$ with
density $f_V$. $P$ obeys the global Markov property with respect to $\G
$ if and only if for every ancestral set $A \in\mathcal{A}(\G)$, and
$\mu$-almost all $x_A \in\X_A$.
%
\begin{eqnarray}
f_A(x_A) &=& \prod_{H\in[A]_{\G}}
f_{H|T}(x_H | x_T). \label{eqnfactorization}
\end{eqnarray}
\end{teo}

A formal proof of this result is given in the \hyperref[app]{Appendix}; a sketch proof
is given in \citet{richardson09}, Theorem~4.

%
\begin{exm} \label{exmfactor1}
For the graph in Figure~\ref{fignoorder}(a), observe that the global
Markov property implies precisely that $X_3 \indep X_4 | X_{12}$, and
$X_1 \indep X_2$. Applying the partition function to the relevant sets
of vertices yields
\begin{eqnarray*}
\bigl[\{1,2,3,4\}\bigr] &=& \bigl\{\{1,4\}, \{2,3\}\bigr\},
\end{eqnarray*}
so Theorem \ref{teofactorization} gives us the factorization from the
\hyperref[sec1]{Introduction}:
\begin{eqnarray*}
f_{1234}(x_1, x_2, x_3,
x_4) &=& f_{23|1}(x_2, x_3 \mid
x_1) \cdot f_{14|2}(x_1, x_4 \mid
x_2)
\end{eqnarray*}
for all $x_i \in\X_i$, $i=1,\ldots,4$.
The expression may appear slightly strange, since the first factor
is the density for $\{X_2, X_3\}$ given $X_1$, while the second is for
$\{X_1, X_4\}$ given $X_2$; nevertheless this factorization does
indeed imply that $X_3 \indep X_4 | X_{12}$. Further, integrating
out $x_3$ and $x_4$ gives
\begin{eqnarray*}
f_{12}(x_1, x_2) &=& f_{2|1}(x_{2}
| x_1) \cdot f_{1|2}(x_{1} | x_2),
\qquad x_1 \in\X_1,  x_2 \in
\X_2,
\end{eqnarray*}
which implies that $X_1 \indep X_2$.
\end{exm}

%
\begin{rmk} \label{rmkmbl}
It follows from Theorem \ref{teofactorization} that if $H$ is a head,
$\tail_\G(H)$ is the Markov blanket for $H$ in the set $\an_\G(H)$, in
the sense that under the global Markov property,
\[
H \indep\an_\G(H) \setminus\bigl(H \cup\tail_\G(H)
\bigr) \mid \tail_\G(H). 
\]
\end{rmk}

%
\begin{rmk}
A different, incorrect definition of $\Phi_\G$ (and, therefore,
$\psi_\G$, $[\cdot]_{\G}$) was given in \citet{richardson09} and
\citet{evans10}. The erroneous definition coincides with that
given here when $W$ is ancestral, so equation~(\ref{eqnfactorization}) holds for both. However, equation
(\ref{eqnpform}) below does not hold for the incorrect partitioning
function in general.
%
\end{rmk}

\section{Toward a parameterization of the discrete Markov model for an ADMG} \label{sectpara}

The factorizations in Theorem \ref{teofactorization} can be used to
produce a parameterization of ADMG models when $\X_V$ is a finite set,
and thus the relevant random variables are discrete. For simplicity
of exposition, we will henceforth assume that the random variables
are\vadjust{\goodbreak}
binary, so $\X_V = \{0,1\}^{|V|}$. Extension to the general finite
discrete case is easy but notationally challenging; this is done in
the special case of ADMGs with chain graph structure by \citet{drton09}.

In the following result, and throughout the paper, empty products are
assumed to equal 1.

%
\begin{teo} \label{teoparametrization}
Let $\G$ be an ADMG, and $P$ a probability distribution on $\{0,1\}
^{|V|}$. Then $P$ obeys the global Markov property with respect to $\G$
if and only if for every ancestral set $A$ and $x_A \in\X_A$,
%
\begin{eqnarray}
P (X_A = x_A ) &=& \sum_{C\dvtx O \subseteq C \subseteq A}
(-1)^{\llvert C \setminus O \rrvert} \prod_{H\in[C]_{\G}} P(X_H =
\vzero| X_T = x_T), \label{eqnpform}
\end{eqnarray}
where $O \equiv\{v \in A | x_v = 0\}$.
\end{teo}

Theorem \ref{teoparametrization} shows that conditional probabilities
of the form $P(X_H = 0 | X_T = x_T)$ are sufficient to form a
parameterization of the binary ADMG model; it remains to show that they
are nonredundant, which is proved in Section~\ref{secsmooth}.

Note that the sets $C$ in (\ref{eqnpform}) may not be ancestral,
which hinders proof by induction. In order to facilitate the proof,
we define the following quantity which will be needed in the
intermediate steps of the induction.

%
\begin{dfn}
Let $A$ be an ancestral set in an ADMG $\G$, and consider a
particular assignment $x_A$ to $X_A$; write $O \equiv\{v \in A
| x_v = 0\}$. For any sets $B$, $W$ such that $B \subseteq W
\subseteq(A\setminus O)$, define the following quantity:
\begin{eqnarray*}
g_{x_A}(B,W) &\equiv &(-1)^{\llvert B \rrvert} \prod
_{H\in[O\cup W]_{\G}} P(X_{H\cap(O\cup B)} = \vzero, X_{H\setminus
(O\cup B)} = \vone\mid
X_T = x_T). 
\end{eqnarray*}
%
\end{dfn}

Note that if $B = \varnothing$ then the right-hand side has factors of
the form $P(X_H = x_H | X_T = x_T)$, and looks much like
(\ref{eqnfactorization}); however, if $B = W$ the expression is a
product of the form $P(X_H = 0 | X_T = x_T)$, just like each term
of (\ref{eqnpform}).

The interpretation is that $W$ is the set of nonzero vertices being
partitioned, and which need to have their values on the left-hand side of
any conditioning bars ``flipped'' from 1 to 0 in order to get an
expression of the form (\ref{eqnpform}). The set $B$ consists of
those vertices for which this ``flipping'' has already taken place, and
those in $W\setminus B$ have yet to be flipped.

The induction starts with the single term $(B,W) = (\varnothing,
A \setminus O)$, given via Theorem \ref{teofactorization}. At each
step a term is ``reduced'' into a sum of two further pieces by flipping
a single vertex, until the procedure finishes with a sum containing
the set of terms
\[
\bigl\{ g_{x_A}(B,W) \dvtx (B,W) = (C,C), \mbox{ where } C \subseteq A
\setminus O \bigr\},
\]
and thus corresponds to an expression of the form (\ref{eqnpform}).

%
\begin{dfn} \label{dfnreduce}
Take a triple $(x_A, B, W)$, where $B \subsetneq W \subseteq(A
\setminus O)$ for $O \equiv\{v \in A | x_v = 0\}$. We say that
$(x_A, B, W)$ is \emph{reducible} if
for each $H \in[O \cup W]_\G$ such that $H \cap(W \setminus B)
\neq\varnothing$, it holds that $\dis_{\an(H)}(H) \setminus H
\subseteq O \cup
(W \setminus B)$.
\end{dfn}

In words, given a set $W$ in which not all vertices are flipped, so
$W\setminus B \neq\varnothing$,
the condition requires that any head $H$ which is in the
partition and has not yet been fully ``flipped,'' has the part of its
tail which is from the same district [i.e., $\dis_{\an(H)}(H) \setminus H$]
consists solely of vertices that are either in $W$ or not yet flipped.






The following technical lemma provides the necessary piece for the
induction step.

%
\begin{lem}\label{lemonestep}
Let $A$ be an ancestral set, and $P$ a distribution obeying the
global Markov property with respect to $\G$. If $(x_A, B, W)$
is reducible in $\G_A$, then there is some $w \in W \setminus B$
such that
%
\begin{equation}
g_{x_A}(B,W) = g_{x_A}\bigl(B, W\setminus\{w\}\bigr) +
g_{x_A}\bigl(B\cup\{w\}, W\bigr), \label{eqngsum}
\end{equation}
and, in addition, either $B = W \setminus\{w\}$ (so $B \cup \{w\} = W$), or both $(x_A, B, W
\setminus\{w\})$ and $(x_A, B \cup\{w\}, W)$ are also reducible.
\end{lem}

\begin{pf}
See the \hyperref[app]{Appendix}.
\end{pf}

Here, $w$ is a vertex that is given the value $1$ in every head in
$g_{x_V}(B, W)$, but is ``flipped'' so it is set equal to $0$ in
$g_{x_V}(B\cup\{w\}, W)$ and is removed from the partition in
$g_{x_V}(B, W\setminus\{w\})$. A major difficulty in the overall
proof of Theorem~\ref{teoparametrization} stems from the fact that,
though each $g_{x_A}$ produced after a reduction is itself reducible or of the form $g_{x_A}(C, C)$, we
will not generally be able to flip the same vertex in each term.

\begin{pf*}{Proof of Theorem \ref{teoparametrization}}
By Theorem \ref{teofactorization}, the global Markov property holds if
and only if for each ancestral $A$ and $x_A$,
\begin{eqnarray*}
P(X_A = x_A) &=& \prod_{H \in[A]_\G}
P(X_H = x_H | X_T = x_T)
\\
&=& g_{x_A}(\varnothing, A \setminus O)
\end{eqnarray*}
using the definition of
$g_{x_V}$. It is easy to check from Definition \ref{dfnreduce} that
either \mbox{$A = O$}, in which case there is nothing to prove, or $(x_A,
\varnothing, A \setminus O)$ is reducible. Then from repeated
application of Lemma
\ref{lemonestep} this is just
\begin{eqnarray*}
&=& \sum
_{C \subseteq A \setminus O} g_{x_A}(C,C)
\\
&=& \sum_{C \subseteq A \setminus O} (-1)^{|C|} \prod
_{H \in[O \cup C]_\G} P(X_H = 0 | X_T =
x_T)
\end{eqnarray*}
which, by inspection, gives the required result.
Conversely, suppose (\ref{eqnpform}) holds, and that $v \in
\barren_\G(A)$ has district $D_1 = \dis_A(v)$;
let $D_2 \equiv A \setminus D_1$ and for $C\subseteq A$, let
$C_i \equiv C \cap D_i$ and \mbox{$O_i \equiv O \cap D_i$},
$i=1,2$.  Then $C \setminus O = (C_1 \setminus O_1)\, \dot\cup\,\break (C_2
\setminus O_2)$, and from Proposition \ref{proppart} get $[C]_\G =
[C_1]_\G \cup [C_2]_\G$.  Hence,
\begin{eqnarray*}
P(X_A = x_A)
&=& \sum_{C\dvtx O \subseteq C \subseteq A} (-1)^{ | C \setminus O
        |}    \prod_{H\in [C]_{\G}} P(X_H = 0 \mid X_T = x_T)
\\
&=& \mathop{\sum_{C_1\dvtx O_1 \subseteq C_1 \subseteq D_1}}_{C_2\dvtx O_2
    \subseteq C_2 \subseteq D_2}       (-1)^{ \mid  C_1 \setminus
     O_1  | +  | C_2 \setminus O_2  |}
\\
&&{}\times     \prod_{H
     \in [C_1]_{\G} \cup [C_2]_\G} P(X_H = 0 \mid X_T = x_T)
     \\
&=& \sum_{C_1\dvtx O_1 \subseteq C_1 \subseteq D_1}       (-1)^{ \mid
      C_1 \setminus O_1  |}    \prod_{H \in [C_1]_{\G}} P(X_H =
      0 \mid X_T = x_T)
\\
&&{} \times\sum_{C_2\dvtx O_2 \subseteq C_2 \subseteq D_2} (-1)^{ |
      C_2 \setminus O_2  |} \prod_{H \in [C_2]_{\G}} P(X_H = 0
      \mid X_T = x_T)
\\
&=& h(x_{D_1}, x_{\pa_\G(D_1)}) \cdot k(x_{D_2}, x_{\pa_\G(D_2)})
\end{eqnarray*}
for some functions $h,k$.  In particular, $k$ does not involve $x_v$,
so it follows that $v \indep A \setminus (D_1 \cup \pa_\G(D_1)) \mid
(D_1 \setminus \{v\}) \cup \pa_\G(D_1)$ which, by the definition of the
Markov blanket of $v$ in $A$, is equivalent to
\[
v \indep A \setminus \bigl(\mbl_\G(v, A) \cup \{v\}\bigr) \mid \mbl_\G(v, A).
\]
It follows that the ordered local Markov property holds (for any
topological ordering); hence, by Proposition \ref{propordered} so
does the global Markov property.
\end{pf*}

\section{Model smoothness}\label{secsmooth}

Let $\mathcal{P}_\G\subseteq\Delta_{2^n-1}$ denote the set of
strictly positive binary probability distributions which obey the
global Markov property with respect to an ADMG $\G$, where
$\Delta_{k}$ is the strictly positive $k$-dimensional probability
simplex and $n$ is the number of vertices in $\G$. We call
$\mathcal{P}_\G$ the \emph{model} defined by $\G$ on a binary
state-space. In this section, such models are shown to be smooth, in
the sense that they are curved exponential families of distributions,
and we prove that the conditional probabilities used in Theorem
\ref{teoparametrization} constitute a parameterization.

Models induced by patterns of conditional independence may be
non-smooth, and determining which are smooth in general is a difficult
open problem [\citet{drton10}]. Non-smoothness can occur even if the
conditional independences arise from a Markov property applied to a
graph, as in the following example.

%
\begin{exm}
Consider the chain graph given in Figure~\ref{figchain}, which mixes
directed and undirected edges. Under the Alternative Markov Property
(AMP) for chain graphs, this graph represents distributions for which
$X_2 \indep X_4 | X_1, X_3$ and $X_1 \indep X_2, X_4$ [\citet
{amp}]. This is shown by \citet{drton09} to represent a non-smooth
model for discrete random variables.

%
\begin{figure}[t]

\includegraphics{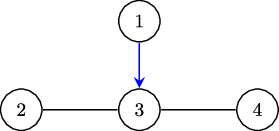}

\caption{A chain graph representing a non-smooth discrete model under
the Alternative Markov Property. (Note that this is not an acyclic
directed mixed graph nor a summary graph.)}\label{figchain}
\end{figure}
\end{exm}

It follows from Theorem \ref{teoparametrization} that for an ADMG $\G$,
the collection of probabilities of the form
\[
P(X_H = \vzero\mid X_T = x_T),\qquad
x_T \in\X_T, H \in\mathcal{H}(\G),
\]
is sufficient to recover the joint distribution under the model
$\mathcal{P}_\G$.
However, it is not immediately clear that each of these probabilities
is necessary, or more specifically that the map in (\ref{eqnpform}) is
smooth and of full rank everywhere.

For brevity, we write $q_H(x_T) \equiv P(X_H = \vzero | X_T = x_T)$,
and the vector of all such probabilities by
%
\begin{equation}
\mathbf{q} \equiv\bigl(q_H(x_T) | H \in
\mathcal{H}(\G), x_T \in\X_T\bigr). \label{eqnq}
\end{equation}
For $\mathbf{p} \in\mathcal{P}_\G$, we---in a mild abuse of notation---let
$\mathbf{q}(\mathbf{p})$ be the vector of the form~(\ref{eqnq})
determined by
calculating the appropriate conditional probabilities from $\mathbf{p}$.
Since this only involves adding and dividing strictly positive
numbers, the map $\mathbf{q}$ is smooth (infinitely differentiable).
Let $\mathcal{Q}_\G\equiv\mathbf{q}(\mathcal{P}_\G)$ be the image of
$\mathbf{q}$ over $\mathcal{P}_\G$;
we call $\mathcal{Q}_\G$ the set of \emph{derived} parameter values. We
will prove that the map in (\ref{eqnpform}) provides a smooth inverse
to $\mathbf{q}$. The first result shows that the set of vectors
${\mathbf{q}}$
that are derived parameters corresponds exactly to those which give
strictly positive probabilities under the inverse map.


%
\begin{teo} \label{teogmpvar}
For an ADMG $\G$, a vector $\mathbf{q}$ is derived (i.e., $\mathbf{q}
\in
\mathcal{Q}_\G$) if and only if for each $x_V \in\X_V$, we have
%
\begin{equation}
p_{x_V}(\mathbf{q}) \equiv\sum_{C\dvtx x_V^{-1} (0) \subseteq C
\subseteq V}
(-1)^{|C \setminus x_V^{-1}(0)|} \prod_{H \in[C]_{\G}} q_{H}(x_T)
> 0, \label{eqnpos}
\end{equation}
where $x_V^{-1}(0) \equiv\{v \in V | x_v = 0\}$.
\end{teo}

%
\begin{rmk} \label{rmkgmpvar}
The boundary of $\mathcal{Q}_\G$ is the set of $\mathbf{q}$ such that
$p_{x_V}(\mathbf{q}) \geq0$ for all $x_V \in\X_V$, with equality
holding in at least one case.

The definition of $p_{x_V}(\mathbf{q})$ in (\ref{eqnpos}) is of the same
form as the expression given for $P(X_V = x_V)$ in (\ref{eqnpform})
and so the result might at first seem trivial; clearly probabilities
must be nonnegative. However, it is not immediately obvious that
this condition is \emph{sufficient} for parameters to be in the
image set $\mathcal{Q}_\G\equiv\mathbf{q}(\mathcal{P}_\G)$.
If we take some $\mathbf{q}^{\dagger} \notin\mathcal{Q}_\G$ and apply to
it the nonlinear functional form in (\ref{eqnpos}) to obtain $\mathbf
{p}(\mathbf{q}^{\dagger})$, without this result there is no apparent reason
why $\mathbf{p}(\mathbf{q}^{\dagger})$ should not be a probability
distribution, nor indeed in
$\mathcal{P}_\G$.
\end{rmk}

To prove Theorem \ref{teogmpvar}, we need the following lemma.

%
\begin{lem} \label{lemsumanc}
Let $A$ be an ancestral set\index{ancestral set} in $\G$, and let $x_A
\in\X_A$. Then for any real vector $\mathbf{q}$ (not necessarily in
$\mathcal{Q}_\G$), the map in (\ref{eqnpos}) satisfies
\begin{eqnarray*}
\sum_{y_V\dvtx y_A = x_A} p_{y_V}(\mathbf{q}) &=& \sum
_{C\dvtx x_A^{-1} (0) \subseteq C \subseteq A} (-1)^{|C \setminus
x_A^{-1}(0)|} \prod
_{H \in[C]_{\G}} q_{H}(x_T),
\end{eqnarray*}
where $x_A^{-1}(0) \equiv\{v \in A | x_v = 0\}$. In particular, taking
$A = \varnothing$,
\[
\sum_{y_V} p_{y_V}(\mathbf{q}) = 1.
\]
\end{lem}

Recall that empty products are assumed equal to 1.

\begin{pf*}{Proof of Lemma \ref{lemsumanc}}
If $A=V$ the result is trivial. If not, pick some $v \in\barren_\G(V)
\setminus A$; this is possible because if $A \supseteq\barren_\G(V)$
then $A = V$ by ancestrality of $A$. So
\begin{eqnarray*}
\mathop{\sum_{y_V\dvtx}}_{y_A = x_A} p_{y_V}(
\mathbf{q}) &=& \mathop{\sum_{y_V\dvtx}}_{y_A = x_A}
\sum_{y_V^{-1} (0) \subseteq C \subseteq V} (-1)^{|C \setminus
y_V^{-1}(0)|} \prod
_{H \in[C]_{\mathcal{G}}} q_{H}(y_T)
\\
&=& \mathop{\sum_{y_{V\setminus\{v\}\dvtx}}}_{y_A = x_A} \sum
_{y_{v}} \sum_{y_V^{-1} (0) \subseteq C \subseteq V}
(-1)^{|C \setminus y_V^{-1}(0)|} \prod_{H \in[C]_{\mathcal{G}}} q_{H}(y_T)
\\
&=& \mathop{\sum_{y_{V \setminus\{v\}\dvtx}}}_{y_A = x_A} \biggl(\sum
_{y_{V \setminus\{v\}}^{-1} (0) \subseteq C \subseteq V} (-1)^{|C
\setminus y_{V \setminus\{v\}}^{-1}(0)|} \prod
_{H \in[C]_{\mathcal{G}}} q_{H}(y_T)
\\
&&\hspace*{15pt}{} + \sum_{y_{V \setminus\{v\}}^{-1}(0) \cup\{
v\} \subseteq C \subseteq V}
(-1)^{|C \setminus(y_{V \setminus\{v\}}^{-1}(0) \cup\{v\})|} \prod_{H
\in[C]_{\mathcal{G}}} q_{H}(y_T)
\biggr).
\end{eqnarray*}
The last equation simply breaks the sum into cases where $y_v = 1$ and
$y_v = 0$, respectively, which takes this form because $v$ does not
appear in any tail sets. The first inner sum in the last expression
can be further divided into the cases where $C$ contains $v$, and
those where it does not, giving
\begin{eqnarray*}
\mathop{\sum_{y_V\dvtx}}_{y_A = x_A} p_{y_V}(
\mathbf{q}) &=& \mathop{\sum_{y_{V \setminus\{v\}}\dvtx}}_{y_A = x_A}
\biggl( \sum_{y_{V \setminus\{v\}}^{-1} (0) \subseteq C \subseteq V
\setminus\{v\}} (-1)^{|C \setminus y_{V \setminus\{v\}}^{-1}(0)|}
\prod
_{H \in[C]_{\mathcal{G}}} q_{H}(y_T)
\\
&&\hspace*{15pt}{}+ \sum_{y_{V \setminus\{v\}}^{-1}(0) \cup\{v\} \subseteq C
\subseteq V} (-1)^{|C \setminus y_{V \setminus\{v\}}^{-1}(0)|} \prod
_{H \in[C]_{\mathcal{G}}} q_{H}(y_T)
\\
&&\hspace*{15pt}{}+ \sum_{y_{V \setminus\{v\}}^{-1}(0) \cup\{v\} \subseteq C
\subseteq V} (-1)^{|C \setminus(y_{V \setminus\{v\}}^{-1}(0) \cup\{v\}
)|} \prod
_{H \in[C]_{\mathcal{G}}} q_{H}(y_T) \biggr).
\end{eqnarray*}
The second and third terms differ only by a factor of $-$1, and so
cancel leaving
\begin{eqnarray*}
\mathop{\sum_{y_V\dvtx}}_{y_A = x_A} p_{y_V}(
\mathbf{q}) &=& \mathop{\sum_{y_{V \setminus\{v\}}\dvtx}}_{y_A = x_A}
\biggl( \sum_{y_{V \setminus\{v\}}^{-1} (0) \subseteq C \subseteq V
\setminus\{v\}} (-1)^{|C \setminus y_{V \setminus\{v\}}^{-1}(0)|}
\prod
_{H \in[C]_{\mathcal{G}}} q_{H}(y_T) \biggr).
\end{eqnarray*}
Repeating this until no vertices outside $A$ are left gives
\begin{eqnarray*}
\mathop{\sum_{y_V\dvtx}}_{y_A = x_A} p_{y_V}(
\mathbf{q}) &=& \sum_{y_{A}^{-1} (0) \subseteq C \subseteq A} (-1)^{|C
\setminus y_{A}^{-1}(0)|}
\prod_{H \in[C]_{\mathcal{G}}} q_{H}(y_T).
\end{eqnarray*}
In the special case $A = \varnothing$, we end up with an empty product
\begin{eqnarray*}
\sum_{y_V} p_{y_V}(\mathbf{q}) &=&
(-1)^{|\varnothing|} \prod_{H \in[\varnothing]_{\mathcal{G}}} q_{H}(y_T)
= 1.
\end{eqnarray*}\upqed
\end{pf*}

\begin{pf*}{Proof of Theorem \ref{teogmpvar}}
The ``only if'' part of the statement follows from Theorem
\ref{teoparametrization} by the fact that if the parameters are
derived then $p_{x_V}(\mathbf{q}) = P(X_V = x_V)$, and these are
therefore positive by definition of $\mathcal{P}_\G$.

For the converse, suppose that the inequalities hold; we will show that
we can retrieve the parameters simply by calculating the appropriate
conditional probabilities. Lemma \ref{lemsumanc} ensures that $\sum
_{x_V} p_{x_V}(\mathbf{q}) = 1$, and that therefore this is a
probability distribution.

Choose some $H^* \in\mathcal{H}(\G)$, with $T^* = \tail_\G(H^*)$ and
$A = \an_\G(H^*)$; also set $x_{H^*} = \vzero$ and pick $x_{T^*} \in\{
0,1\}^{|T^*|}$. By Lemma \ref{lemsumanc},
\begin{eqnarray*}
\sum_{y_V\dvtx y_{A} = x_{A}} p_{y_V}(\mathbf{q}) &=& \sum
_{y_{A}^{-1} (0) \subseteq C \subseteq A} (-1)^{|C \setminus
y_{A}^{-1}(0)|} \prod
_{H \in[C]_{\mathcal{G}}} q_{H}(y_T).
\end{eqnarray*}
Now clearly $H^* \in\Phi_\G(A)$, so applying Lemma
\ref{lemheadfirst} and the fact that $H^* \subseteq x_{A}^{-1}(0) =
y_{A}^{-1}(0)$ shows $H^* \in[C]_\G$ for all terms $C$ in the sum
and, therefore, we can apply Proposition \ref{proprmv} to factor out
the parameter associated with~$H^*$:
\begin{eqnarray*}
&=& q_{H^*}(y_{T^*}) \sum_{y_{A}^{-1} (0) \subseteq C \subseteq A}
(-1)^{|C \setminus y_{A}^{-1}(0)|} \prod_{H \in[C \setminus
H^*]_{\mathcal{G}}} q_{H}(y_T)
\\
&=& q_{H^*}(y_{T^*}) \sum_{y_{A \setminus H^*}^{-1} (0) \subseteq C
\subseteq A \setminus H^*}
(-1)^{|C \setminus y_{A \setminus H^*}^{-1}(0)|} \prod_{H \in
[C]_{\mathcal{G}}} q_{H}(y_T).
\end{eqnarray*}
But note that $A \setminus H^*$ is also an ancestral set, and thus
using Lemma \ref{lemsumanc} again,
\begin{eqnarray*}
\sum_{y_V\dvtx y_{A \setminus H^*} = x_{A \setminus H^*}}
p_{y_V}(\mathbf{q}) &=& \sum
_{y_{A \setminus H^*}^{-1} (0) \subseteq C \subseteq A \setminus H^*}
(-1)^{|C \setminus y_{A \setminus H^*}^{-1}(0)|} \prod
_{H \in[C]_\mathcal{G}} q_{H}(y_{T}).
\end{eqnarray*}
Hence,
\begin{eqnarray*}
\frac{\sum_{y_{V \setminus A}} p_{x_V}(\mathbf{q})}{\sum_{y_{V
\setminus(A \setminus H^*)}} p_{x_V}(\mathbf{q})} &=& q_{H^*}(x_{T^*}),
\end{eqnarray*}
and we can recover the original parameters from the probability
distribution $\mathbf{p}$ in the manner we would expect; that $\mathbf{p}$
satisfies the global Markov property for $\G$ then follows from
Theorem \ref{teoparametrization}. Thus, $\mathbf{p} \in\mathcal{P}_\G$
and $\mathbf{q} = \mathbf{q}(\mathbf{p}) \in\mathcal{Q}_\G$.
\end{pf*}

%
\begin{teo}
For an ADMG $\G$, the model $\mathcal{P}_\G$ of strictly positive
binary probability distributions satisfying the global Markov property
with respect to $\G$ is smoothly parameterized by $\mathbf{q} \in
\mathcal{Q}_\G$.

Consequently, the model $\mathcal{P}_\G$ is a curved exponential family
of dimension
\[
d= \sum_{H \in\mathcal{H}(\G)} |\X_{\tail(H)}| =
\sum_{H \in\mathcal{H}(\G)} 2^{| \tail(H)|}.
\]
\end{teo}

\begin{pf}
By Theorem \ref{teogmpvar}, the set $\mathcal{Q}_\G\subseteq\reals^d$
is open. The map $\mathbf{p}(\mathbf{q})\dvtx \mathcal{Q}_\G
\rightarrow\mathcal{P}_\G$ is multilinear and, therefore, infinitely
differentiable. Its inverse $\mathbf{q}\dvtx\mathcal{P}_\G\rightarrow
\mathcal{Q}_\G$ is also infinitely differentiable on $\mathcal{P}_\G$.

The composition $\mathbf{q} \circ\mathbf{p}$ is the identity function
on $\mathcal{Q}_\G$ and, therefore, its Jacobian is the identity matrix
$I_d$. However, the Jacobian of a composition of differentiable
functions is the product of the Jacobians, so
\[
I_d = \frac{\partial\mathbf{q}}{\partial\mathbf{p}} \frac{\partial
\mathbf{p}}{\partial\mathbf{q}}.
\]
But this implies that each of the Jacobians has full rank $d$ and,
therefore, the map $\mathbf{q}$ is a smooth parameterization of
$\mathcal{P}_\G$. See \citet{kass97}, Corollary A.3.
\end{pf}

\section{Discussion}\label{secdiscussion}
We remark that it is easy to extend the results of Sections~\ref
{sectpara} and \ref{secsmooth} from the binary case to a general finite
discrete state-space; we have avoided this only for notational simplicity.
It is also a simple matter to extend the results from ADMGs to the \emph
{summary graphs} of \citet{wermuth11} which incorporate three
types of edge:
directed ($\rightarrow$),
undirected ({\tikz[baseline=-3pt] \draw(0cm,0cm) -- (5mm,0cm);}),
and dashed ({\tikz[baseline=-3pt] \draw[dashed] (0cm,0cm) --
(5mm,0cm);}); the dashed edges are equivalent to bidirected
($\leftrightarrow$) edges [\citet{sadeghi11}].
The undirected component of a summary graph can be dealt with using
standard methods for undirected graphs [\citet{lau96}], and the
remaining parameterization done as for an ADMG, conditional on the undirected
component.

\begin{appendix}\label{app}
\section*{Appendix: Technical proofs}\vspace*{-6pt}\label{appA}
\begin{pf*}{Proof of Lemma \ref{lempartsuit}}
Suppose that two heads $H_1, H_2$ are distinct and \mbox{$H_1 \cap H_2
\neq\varnothing$}. We will show that they are dominated by $H^*
\equiv\barren_\G(H_1 \cup H_2)$; clearly $H^* \subseteq H_1 \cup
H_2$ and $H_1, H_2 \subseteq\an_\G(H^*)$, so if $H^*$ is a head
then $\prec$ satisfies the requirements for partition-suitability.

Clearly $H^*$ is barren, so we need to prove that it is contained
within a single district in $\an_\G(H^*)$. By definition,
$\an_\G(H^*) \supseteq H_1 \cup H_2$; we need to find a bidirected
path between any $v, w \in H^* \subseteq H_1 \cup H_2$. If $v$ and
$w$ are either both in $H_1$ or both in $H_2$, then the existence of
such a path follows from the fact that these are heads. If $v \in
H_1$ and $w \in H_2$, then construct a bidirected path in
$\an_\G(H_1)$ from $v$ to some vertex $x \in H_1 \cap H_2$, and a
bidirected path in $\an_\G(H_2)$ from $x$ to $w$; these paths can
then be concatenated into a new path meeting the requirements,
shortening the resulting sequence of edges if necessary to avoid
repetition of vertices. Hence, $H^*$~is a head.

Since $H_1, H_2 \subseteq\an_\G(H^*)$ we have $H_i \preceq H^*$ for
each $i=1,2$, and therefore $\prec$ is partition-suitable.
\end{pf*}

\subsection*{Proof of factorization}

%
\begin{prop} \label{propequiv}
Let $\prec$ and $\prec'$ be two partition-suitable partial orderings
for $\mathcal{H}$, such that\vspace*{1pt} for every $H\in\mathcal{H}$ and $W
\subseteq V$, $H$ is maximal in $W$ under $\prec$ whenever this is so
under $\prec'$. Then $[\cdot]^{\prec} = [\cdot]^{\prec'}$.
\end{prop}

\begin{pf}
We again proceed by induction on the size of $W$. Recall that for
all $v \in V$, we have $\{v\} \in\mathcal{H}$ by the definition of
partition-suitability, so $[\{v\}]^{\prec} = [\{v\}]^{\prec'} =
\{\{v\}\}$. Now take a general $W \subseteq V$, and suppose that
$H$ is maximal under $\prec'$ in $W$; then by Proposition
\ref{proprmv}
\begin{eqnarray*}
[W]^{\prec'} &=& \{H\} \cup[W \setminus H]^{\prec'}
\\
&=& \{H\} \cup[W \setminus H]^{\prec}
\\
&=& [W]^{\prec}
\end{eqnarray*}
by applying the induction hypothesis to $W \setminus H$, and using the
fact that $H \in[W]^\prec$ because it is also maximal under $\prec$
in $W$.
\end{pf}

Define a partial ordering $\prec^*$ on heads in an ADMG by $H_1
\prec^* H_2$ if and only if both $H_1 \prec H_2$, and $H_1$ and $H_2$
are contained in the same district in $\an_\G(H_1 \cup H_2)$; note
that this is a weaker ordering than $\prec$, since strictly fewer
pairs of sets are comparable. It is easy to see that $\prec^*$ is
partition-suitable for heads $\mathcal{H}(\G)$ by repeating the proof
of Lemma \ref{lempartsuit}. In addition, sets which are maximal
under $\prec$ will also be maximal under $\prec^*$, so the partitions
defined by these two orderings are the same by Proposition
\ref{propequiv}.

This weaker partial ordering leads us to a class of sets which play
a role similar to that of ancestral set: a set with ``ancestrally
closed districts'' is one whose districts are ancestrally closed
(rather than the whole set).

%
\begin{dfn} \label{dfnacd}
Let $\G$ be an ADMG, and $W$ be a subset of its vertices. We say
$W$ \emph{has ancestrally closed districts} if $\dis_{\an(W)}(W) =
W$.
\end{dfn}

Equivalently, $W$ has ancestrally closed districts if $W$ is not
connected to $\an_\G(W) \setminus W$ by any bidirected edges. This
definition is important because the partitioning function $[\cdot]_\G$
will act upon sets with ancestrally closed districts ``separately''
within the relevant ancestral set: that is, for such sets,
\[
\bigl[\an_\G(W)\bigr]_\G= [W]_\G\cup
\bigl[\an_\G(W) \setminus W\bigr]_\G.
\]
Note that if \mbox{$D = D_1\,\dot\cup\, D_2$} has ancestrally closed
districts, and $D_1$ and $D_2$ are not joined by any bidirected edges,
then $D_1$ and $D_2$ themselves have ancestrally closed districts
(here $\dot{\cup}$ indicates a disjoint union). If for every $v,w \in
D$ there is a bidirected path from $v$ to $w$ such that all the
vertices on the path are contained within $D$, then $D$ cannot be
partitioned in this manner, and we say it is
\emph{bidirected-connected}.


%
\begin{dfn}
Let $C \subseteq V$. We say that an ordering $<$ on the vertices of
$C$ is \emph{$(C, \prec^*)$-consistent} if for any $H_1, H_2 \in
[C]_\G$ such that $H_1 \prec^* H_2$, we have $v_1 < v_2$ for all
$v_1 \in H_1$, $v_2 \in H_2$.
\end{dfn}

%
\begin{lem} \label{lemdisjoint-extend} Let \mbox{$D = D_1\,\dot\cup\, D_2$}
have ancestrally closed districts and be such that $D_1$ is not
connected to $D_2$ by any bidirected edges. Let $<_1$ and $<_2$ be
orderings on $D_1$ and $D_2$ (resp.). If for $i=1,2$, $<_i$
is $(D_i,\prec^*)$-consistent, then every extension of $<_1$ and
$<_2$ to an ordering $<$ on $D$ is also a $(D, \prec^*)$-consistent
ordering.
\end{lem}

\begin{pf}
Orderings between vertices $v_1,v_2 \in D_i$ are specified by
$<_i$. Further, if $v_1 \in D_1$ and $v_2 \in D_2$ then since $v_1$
and $v_2$ are in different districts in $\an_\G(D)$, it follows from
the definition of $\prec^*$ that $v_1$ and $v_2$ can be ordered in
any way to achieve a $(D, \prec^*)$ consistent ordering.
\end{pf}

A total ordering $<_i$ on a set $D_i$ will be said to be {\em
topological in $\G$} if no vertex $d\in D_i $ precedes any of its
proper ancestors in $\G$ that are in $D_i$.

%
\begin{lem} \label{lemdisjoint-extend-top}
Let $D_1$ and $D_2$ be disjoint subsets in $\G$. Let $<_1$ and
$<_2$ be topological orderings on $D_1$ and $D_2$ (resp.).
Then there exists an extension of $<_1$ and $<_2$ to a topological
ordering $<$ on $D_1 \cup D_2$.
\end{lem}

\begin{pf}
We construct a topological ordering iteratively as follows: let
$\langle d_1,\ldots, d_{k-1}\rangle$ be the first $k-1$ vertices in
$D_1 \cup D_2$ already ordered under $<$; let $E_k = (D_1 \cup D_2)
\setminus\{d_1,\ldots, d_{k-1}\}$ be the set of vertices remaining
to be ordered. Further, let $Q_k = \{ d \mid d \in E_k, \an_{\G}(d)
\cap E_k = \{d\}\}$ be those vertices in $E_k$ that have no proper
ancestors in $E_k$; $Q_k \neq\varnothing$ since $V$ is finite and
$\G$ is acyclic. Finally, if $Q_k \cap D_1 \neq\varnothing$, define
$d_k$ to be the first element in $Q_k$ under $<_1$, otherwise define
$d_k$ to be the first element in $Q_k$ under $<_2$. That the
ordering is topological follows from the definition of $Q_k$.
\end{pf}

%
\begin{lem} \label{lemancdist}
Let $D$ have ancestrally closed districts, and suppose $C\subseteq
\barren_\G(D)$.
Then $D \setminus C$ has ancestrally closed districts.
\end{lem}

\begin{pf}
Let $D' \equiv D \setminus C$. Since $C\subseteq\barren_\G(D)$, $\an_\G
(D \setminus C) \subseteq\an_\G(D) \setminus C$,
so
\begin{eqnarray*}
\dis_{\an(D')}\bigl(D'\bigr) \subseteq\dis_{\an(D)}
\bigl(D'\bigr) \setminus C \subseteq\dis_{\an(D)}(D)
\setminus C = D \setminus C &=& D'.
\end{eqnarray*}
Since $D' \subseteq\dis_{\an(D')}(D')$, the result holds.
\end{pf}

%
\begin{lem}\label{lemcw-extend}
Let $C \cup W$ have ancestrally closed districts,
with $W \subseteq\barren_\G(C \cup W)$ and $W \cap C = \varnothing$.
Then any ordering on $W$ may be extended to a topological ordering
of the vertices in $C \cup W$ which is both $(C, \prec^*)$ and $(C
\cup W, \prec^*)$-consistent.
\end{lem}

\begin{pf}
Note that $C$ has ancestrally closed districts by Lemma
\ref{lemancdist}. We proceed by induction on the size of
$C \cup W$; if $|C \cup W| = 0$ or $1$ then the result is trivial.

If $C \cup W$ contains two components which are not connected by
bidirected edges, then we can split it into two smaller sets $C_1
\cup W_1$ and $C_2 \cup W_2$, each with ancestrally closed districts, where $C = C_1\,\dot{\cup}\, C_2$ and $W
= W_1\,\dot{\cup}\, W_2$.
Clearly, $W_i \in
\barren(C_i \cup W_i)$ for each $i$, so using the induction
hypothesis, we can find topological orderings $<_i$ on the vertices
of $C_i \cup W_i$ which are both $(C_i \cup W_i, \prec^*)$ and
$(C_i, \prec^*)$ consistent. It then follows from Lemma
\ref{lemdisjoint-extend-top}, taking $D_i = (C_i \cup W_i)$, that
there exists a topological ordering $<$ on $C \cup W$ that extends
$<_1$ and $<_2$. It further follows from two applications of Lemma
\ref{lemdisjoint-extend} that $<$ is both $(C, \prec^*)$ and $(C
\cup W, \prec^*)$-consistent.

Since, by assumption,  $C \cup W$ has ancestrally closed districts, if
this set does not contain two components then $C \cup W$ is a single
district in\break  $\an_\G(C \cup W)$. Let $H = \barren_\G(C \cup W)$;
this is clearly a head and maximal under $\prec^*$ in $C \cup W$.
Further, $W \subseteq H$ so applying Proposition \ref{proprmv} gives
\begin{eqnarray*}
[C \cup W]_\G&=& \{H\} \cup\bigl[(C \cup W) \setminus H
\bigr]_\G
\\
&=& \{H\} \cup\bigl[C \setminus(H \setminus W)\bigr]_\G,
\end{eqnarray*}
since $W \cap C = \varnothing$. Since $H \setminus W \subseteq
\barren_\G(C)$, Lemma \ref{lemancdist} shows that
$C \setminus\break (H \setminus W)$ also has
ancestrally closed districts; applying the induction hypothesis, we
can find a topological ordering of $C$ which is both $(C \setminus(H
\setminus W), \prec^*)$ and \mbox{$(C, \prec^*)$-}consistent [possibly $C
\setminus(H \setminus W) = C$ in which case this is trivial]. This
ordering may be combined with an arbitrary ordering on $W$ by simply
concatenating the orderings so that everything in $W$ comes after
everything in $C$. This gives an ordering which is $(C \cup W,
\prec^*)$-consistent, because $H \supseteq W$ is maximal; since $W$ is
barren in $C \cup W$, the ordering is also topological.
\end{pf}

%
\begin{cor}\label{coradd-one}
If $D \cup\{w\}$ has ancestrally closed districts with $w \in\barren
_\G(D \cup\{w\})$, then there exists an ordering $<$ which is both
$(D, \prec^*)$ and $(D \cup\{w\}, \prec^*)$-consistent, and such that
$w$ is the maximal vertex under $<$.
\end{cor}

\begin{pf}
The claim is trivial if $w \in D$. Otherwise, $\{w\}$ is barren in
$D \cup\{w\}$, so we apply the previous lemma. 
\end{pf}

Note that the previous lemma and this corollary do not generalize to
adding two vertices: there exist graphs with ancestral sets $A$,
$A\cup\{w_1\}$ and $A\cup\{w_1, w_2\}$, such~that no topological
ordering is $(A, \prec^*)$-, $(A \cup\{w_1\}, \prec^*)$- and $(A \cup
\{w_1,w_2\}, \prec^*)$-consistent. See \citet{richardson09} for such
an example.

Given a path, $\pi$, and two vertices $v,w$ on $\pi$, the
\emph{subpath} $\pi(v,w)$ is the sequence of edges which lie between
$v$ and $w$ on $\pi$. As with a path, we allow a single vertex (and
no edges) to be a degenerate case of a subpath.

%
\begin{lem} \label{lempaths}
Suppose $\pi$ is a path from $a$ to $b$, and is not blocked by $C$.
Then every vertex $v$ on $\pi$ is contained in $\an_{\G}(\{a,b\}
\cup C)$.
\end{lem}

\begin{pf}
Suppose $w$ is on ${\pi}$ and is an ancestor of neither $a$ nor
$b$. Then on each of the subpaths ${\pi}( a,w)$ and ${\pi}(w, b)$,
there is at least one edge with an arrowhead pointing towards $w$
along the subpath. Let $v_{aw}$ and $v_{wb}$ be the vertices at
which such arrowheads occur that are closest to $w$ on the
respective subpaths.
There are now three cases:
(1)~if $w \neq v_{wb}$ then ${\pi}(w,v_{wb})$ is a directed path
from $w$ to $v_{wb}$. It further follows that $v_{wb}$ is a collider
on ${\pi}$, and since the path is not blocked by $v_{wb}$, it is an
ancestor of $C$. Hence, $w \in\an_\G(C)$.
(2) If $w \neq v_{aw}$, then a symmetric argument holds.
(3) If $v_{aw} = w = v_{wb}$, then
$w$ is a collider on ${\pi}$, hence again an ancestor of $C$.
\end{pf}

The next two lemmas are used to establish necessary results about
Markov blankets:

%
\begin{lem} \label{lemords} 
Let $H_1,H_2 \in[D]_\G$ with $H_1 \neq H_2$. Then at
least one of the following holds:
\begin{longlist}[(iii)]
\item[(i)] $H_1 \prec H_2$;

\item[(ii)] $H_2 \prec H_1$; or

\item[(iii)] there is no bidirected path between any $h_1 \in H_1$ and $h_2
\in H_2$ contained within $\an_\G(H_1 \cup H_2)$.
\end{longlist}
\end{lem}

\begin{pf}
Suppose $H_1,H_2 \in[D]_\G$, and that (iii) fails. Then let $H^*
\equiv\barren_\G(H_1 \cup H_2)$. Since, $H_1, H_2$ are heads and
since (iii) fails, $H^*$ is a barren set which is connected by
bidirected paths in $\an_\G(H^*) = \an_\G(H_1 \cup H_2)$; hence, it
is a head. In addition, $H^* \subseteq H_1 \cup H_2 \subseteq D$,
and $H^* \succeq H_1, H_2$.

It follows that $H^* \in[D]_\G$, which means that either $H^* =
H_1$, in which case (ii) holds, or $H^* = H_2$, in which case (i)
holds.
\end{pf}



%
\begin{lem} \label{lemmb2}
Let $D$ be bidirected-connected with ancestrally closed
districts, and let $D' \equiv D \setminus\{w\}$ for some $w \in
\barren_\G(D)$. Let $<$ be a total order that is $(D, \prec^*)$- and
$(D', \prec^*)$-consistent, and under which $w$ is maximal. For
a given $v \in D'$ define $H,H'$ to be the heads such that
$v \in H \in[D]_{\G}$ and $v \in H' \in[D']_{\G}$, respectively, and
$T,T'$ the corresponding tails.
Let
\begin{eqnarray*}
B &\equiv& \bigl(\dis_{\pre_<(v)}(v) \setminus\{v\}\bigr) \cup
\pa_\G\bigl(\dis_{\pre_<(v)}(v)\bigr),
\\
C &\equiv& \bigl(H \cap\pre_<(v)\bigr) \cup T
\quad\mbox{and}
\\
C' &\equiv& \bigl(H' \cap\pre_<(v)
\bigr) \cup T'.
\end{eqnarray*}
Then $B \subseteq C$ and $B \subseteq C'$, and $B$ \mbox{m-}separates $v$ from
both $C \setminus B$ and $C' \setminus B$.
%
\end{lem}



\begin{pf}
%
Let $S \equiv\dis_{\pre_<(v)}(v) \subseteq D$; we claim that $S
\subseteq\an_\G(H)$. If not then there is a bidirected path $\pi$
from $v$ to some $s \in S \setminus\an_\G(H)$; let this path be
minimal, so that $s$ is adjacent on $\pi$ to some $t \in\an_\G(H)$.
Then $s$ lies in some different head $H^* \in[D]_\G$, and we have
constructed a bidirected path from $H$ to $H^*$ within $\an_\G(H
\cup H^*)$; it follows from Lemma \ref{lemords} that either $H
\prec H^*$ or $H^* \prec H$, but the former is ruled out by the
existence of $\pi$ and the $(D, \prec^*)$-consistency of $<$. Hence,
$H^* \subseteq\an_\G(H)$, and in particular $s \in\an_\G(H)$, so
we reach a contradiction.\vadjust{\goodbreak}

Thus, $S \subseteq\an_\G(H)$ and, therefore, $S \equiv
\dis_{\pre_<(v)}(v) \subseteq\dis_{\an(H)}(v)$, so
\[
S \cup\pa_\G(S) \subseteq\dis_{\an(H)}(v) \cup\pa\bigl(
\dis_{\an(H)}(v)\bigr) \subseteq H \cup T.
\]
%
%
%
Finally, using $S \subseteq\pre_<(v)$, we have
\[
B = \bigl(S \setminus\{v\}\bigr) \cup\pa_\G(S) \subseteq\bigl(H
\cap\pre_<(v)\bigr) \cup T=C.
\]
It follows from Lemma \ref{lemancdist} and the fact that $w \in
\barren_\G(D)$, that $D'$ also has ancestrally closed districts, and
the same argument as above shows that $B \subseteq C'$.

Now, let $\pi$ be a path from $v$ to some $c \in C \setminus B$, and
assume without loss of generality that $\pi$ does not intersect $C
\setminus B$ other than at $c$. We will show that $\pi$ is blocked by
$B$.

Note that $B \subseteq C \subseteq\pre_{<}(v)$; thus if $\pi$
includes any vertex $s > v$ then it is blocked by Lemma
\ref{lempaths}, because $s$ is not an ancestor of any element of $C$.
Consequently, we may assume that the edge on $\pi$ adjacent to $v$ is
of the form $v \leftrightarrow$ or $v \leftarrow$.

We claim that $\pi$ contains at least one non-collider; suppose not
for a contradiction: then $\pi$ is of the form
\[
v \leftrightarrow t_1 \leftrightarrow\cdots\leftrightarrow
t_p \leftrightarrow c, \qquad v \leftrightarrow t_1
\leftrightarrow\cdots\leftrightarrow t_p \leftarrow c\quad\mbox{or}\quad v \leftarrow c,
\]
with every node $t_i$ an ancestor of $B$ and hence of $D$. Since $D$
has ancestrally closed districts, it follows that every $t_i \in D$
and hence\vspace*{1pt} $t_i \in\dis_{\pre_<(v)}(v) \setminus\{v\}$, so $t_i \in
B$. But then $c \in B$, which is a contradiction, since we
assumed $c \in C \setminus B$.

It follows that $\pi$ contains at least one non-collider; let $d$ be
the non-collider closest to $v$ on the path. But then repeating the
argument above (replacing $c$ with~$d$) shows that $d \in B$ and,
therefore, $\pi$ is blocked by $B$.

Similarly, all paths $\pi'$ from $v$ to some $c'$ in $C' \setminus
B$ are blocked by $B$.
\end{pf}

The next lemma is the crux of the induction used in the proof of
Theorem~\ref{teofactorization}.

%
\begin{lem} \label{lemrmv-one}
Let $D$ have ancestrally closed districts, and $w \in\break 
\barren_\G(D)$. Then for any $f_V$ obeying the global Markov
property with respect to $\G$, we have
\[
\prod_{H \in[D]_\G} f_{H|T}(x_H |
x_T) = f_{w|\an(D) \setminus\{w\}}(x_w | x_{\an(D) \setminus\{w\}})
\prod
_{H \in[D \setminus\{w\}]_\G} f_{H|T}(x_H |
x_T)
\]
$\mu$-almost everywhere.
\end{lem}

\begin{pf}
Note that we need only prove the case where $D$ forms a single
district, from which the general result will follow because by
Proposition \ref{proppart} the factors not involving $\dis_D(w)$
are the same on both sides. Assume therefore that $D = \dis_D(w)$,
and thus $D$ is bidirected-connected.

Define $D' = D \setminus\{w\}$, and let $<$ be a topological total
ordering which is $(D, \prec^*)$ and $(D', \prec^*)$ consistent,
which exists by Corollary \ref{coradd-one}. Further, we can choose
$w$ to be the maximal element in $D$.

For any $v \in H \in[D]_\G$, let $H_v = H \cap\pre_<(v)$, and
similarly for $v \in H' \in[D']_\G$, let $H_v' = H' \cap\pre_<(v)$.
In addition, let
\[
B_v \equiv\bigl(\dis_{\pre_<(v)}(v) \setminus\{v\}\bigr) \cup
\pa_\G\bigl(\dis_{\pre_<(v)}(v)\bigr).
\]
Then
\begin{eqnarray*}
\prod_{H \in[D]_\G} f_{H|T}(x_H |
x_T)
&=& \prod_{H \in[D]_\G} \prod
_{v \in H} f_{v | H_v \cup T}(x_v | x_{H_v},
x_T)
\\
&=& \prod_{H \in[D]_\G} \prod_{v \in H}
f_{v|B_v}(x_v | x_{B_v})
\\
&=& \prod_{v \in D} f_{v|B_v}(x_v |
x_{B_v}),
\end{eqnarray*}
where the first equality follows from the elementary properties of
conditional probabilities, and the second from applying Lemma \ref
{lemmb2} to see that $B_v$ \mbox{m-}separates $v$ from $(H_v \cup T)
\setminus B_v$.

But $B_v$ also \mbox{m-}separates $v$ from $(H'_v \cup T') \setminus B_v$, so
reversing the argument gives
\begin{eqnarray*}
\prod_{v \in D} f_{v|B_v}(x_v | x_{B_v})
&=& f_{w|B_w}(x_w | x_{B_w}) \prod
_{v \in D \setminus\{w\}} f_{v|B_v}(x_v | x_{B_v})
\\
&=& f_{w|B_w}(x_w | x_{B_w}) \prod
_{H' \in[D \setminus\{w\}]_\G} \prod_{v \in H'}
f_{v|B_v}(x_v | x_{B_v})
\\
&=& f_{w|B_w}(x_w | x_{B_w}) \prod
_{H' \in[D \setminus\{w\}]_\G} \prod_{v \in H'}
f_{v|H_v' \cup T'}(x_v | x_{H_v'}, x_{T'})
\\
&=& f_{w|B_w}(x_w | x_{B_w}) \prod
_{H' \in[D \setminus\{w\}]_\G} f_{H' | T'}(x_{H'} |
x_{T'}).
\end{eqnarray*}
In addition, note that $B_w = H_w \cup T$, so it is the Markov blanket
for $w$ in $\an_\G(D)$ using the ordered local Markov property. Thus,
\[
f_{w|B_w}(x_w | x_{B_w}) = f_{w|\an(D) \setminus\{w\}}(x_w
| x_{\an(D) \setminus\{w\}}),
\]
which gives the result.
\end{pf}

\begin{pf*}{Proof of Theorem \ref{teofactorization}}
We proceed by induction on $|A|$. Clearly, the result holds if $|A|
=1$.

If $|A| > 1$, then let $w \in\barren_\G(A)$; thus $A' \equiv A
\setminus\{w\}$ is also ancestral. Suppose that the global Markov
property holds; then by elementary laws of probability and the
induction hypothesis,
\begin{eqnarray*}
f_A(x_A) &=& f_{w | A'}(x_w |
x_{A'}) \cdot f_{A'}(x_{A'})
\\
&=& f_{w | A'}(x_w | x_{A'}) \prod
_{H' \in[A']_\G} f_{H'|T'}(x_{H'} |
x_{T'})
\end{eqnarray*}
and by Lemma \ref{lemrmv-one}, this is just
\[
\hspace*{-38pt}= \prod
_{H \in[A]_\G} f_{H|T}(x_{H} |
x_{T}).
\]
Conversely, suppose that (\ref{eqnfactorization}) holds and let $<$
be a topological ordering of the ancestral set $A$. By the induction
hypothesis, the ordered local Markov property is satisfied for $<$ and
all suitable pairs $(v, A')$ such that $A' \subset A$.
Let $w \in\barren_\G(A)$ be the maximal vertex under $<$ in $A$,
with $H$ such that $w \in H \in[A]_\G$;
the factorization implies that $H \indep A \setminus(H \cup T) \mid T$.
Note that $H = \barren_\G(\dis_A(w))$,~so
\begin{eqnarray*}
\mbl_\G(w,A) &\equiv&\bigl(\dis_A(w) \setminus\{w\}
\bigr) \cup\pa_\G\bigl(\dis_A(w)\bigr)
\\
&=& \bigl(H \setminus\{w\}\bigr) \cup T.
\end{eqnarray*}
This then implies
$w \indep A \setminus(\mbl_\G(w) \cup\{ w\}) | \mbl_\G(w)$
by the weak union property of conditional independence.
Hence, the ordered local Markov property is satisfied.
\end{pf*}

\subsection*{Proof of parameterization}

%
\begin{prop}\label{propanc-dist-part}
If $H \in[W]_{\G}$ and $ D = \dis_{\an(H)}(H) \cap W$ then
$[W]_{\G} = [W\setminus D]_{\G} \cup[D]_{\G}$.
\end{prop}

\begin{pf} Note that since $H \in[W]_{\G}$, $H\subseteq \dis_{\an
(H)}(H) \cap W = D$.
The proof is by induction on $|W\setminus D|$. If $W\setminus
D=\varnothing$, the claim is trivial.
Suppose $H^* \in[W]_{\G}$ and $H^* \cap D \neq\varnothing$.
Applying Lemma \ref{lemords} to $H, H^* \in[W]_\G$ we see that either
$H^* = H$ or $H^* \prec H$,
so $H^* \subseteq D$. Thus, every head in $[W]_\G$ is either a subset
of $D$ or $W\setminus D$.
Consequently,
there exists $H^\dagger\in[W]_\G$ with $H^\dagger\subseteq
W\setminus D$; let $W^\dagger\equiv W\setminus H^\dagger$.
By Proposition \ref{proprmv}, $[W]_{\G} = \{H^\dagger\} \cup[W^\dagger
]_{\G}$.
Since $D\subseteq W^\dagger$ and $H \in[W^\dagger]_{\G}$, the conclusion
follows from the inductive hypothesis applied to $W^\dagger$.
\end{pf}



\begin{pf*}{Proof of Lemma \ref{lemonestep}}
It suffices to prove the result for $A = V$, from which the general
case follows by applying it to the subgraphs $\G_A$.

Since $(x_V, B, W)$ is reducible, $W \setminus B \neq\varnothing$;
let $H^*$ be a maximal head such that both $H^* \in[O \cup W]_\G$
and $H^* \cap(W \setminus B) \neq\varnothing$, further take $w \in
H^* \cap
(W \setminus B)$. Let $D^* \equiv\dis_{\an(H^*)}(H^*)$ be the
associated district within the ancestors of $H^*$. By
construction, $D^*$ has
ancestrally closed districts and is bidirected-connected.

Define $y_B \equiv0$, $y_{V \setminus B} \equiv x_{V \setminus B}$; then
\begin{eqnarray*}
\hspace*{-5pt}&& g_{x_V}(B,W)
\\
\hspace*{-5pt}&&\qquad \equiv (-1)^{\llvert B \rrvert} \prod
_{H\in[O\cup W]_{\G}}\! P(X_{H} = y_H | X_T
= x_T)
\\
\hspace*{-5pt}&&\qquad  = (-1)^{\llvert B \rrvert} \prod_{H\in[(O\cup W) \setminus H^*]_{\G}}\!
P(X_{H} = y_H | X_T = x_T)
\\
\hspace*{-5pt}&&\hspace*{132pt}{}\times\bigl\{ P(X_w = 1, X_{{H^*}\setminus\{w\}} = y_{{H^*}\setminus
\{w\}} |
X_{T^*} = x_{T^*})
\\
\hspace*{-5pt}&&\hspace*{148pt}{}+ P(X_w = 0, X_{{H^*}\setminus\{w\}} = y_{{H^*}\setminus\{w\}} |
X_{T^*} = x_{T^*})
\\
\hspace*{-5pt}&&\hspace*{148pt}{} - P(X_w = 0, X_{{H^*}\setminus\{w\}} = y_{{H^*}\setminus\{w\}} |
X_{T^*} = x_{T^*}) \bigr\}.
\end{eqnarray*}
The last term after distributing the product is just $g_{x_V}(B \cup
\{w\}, W)$, so to prove (\ref{eqngsum}) we need to show that
%
\begin{eqnarray}\label{eqnsuffice1}
&& g_{x_V}\bigl(B, W\setminus\{w\}\bigr)\nonumber
\\[3pt]
&&\qquad \equiv  (-1)^{|B|}
\prod_{H \in[(O \cup W)\setminus\{w\}]_\G}\! P(X_H = y_H |
X_T = x_T)\nonumber
\\[3pt]
&&\qquad = (-1)^{\llvert B \rrvert}
\nonumber\\[-7pt]\\[-7pt]
&&\quad\qquad{}\times \prod_{H \in[(O\cup W) \setminus H^*]_{\G}}\!
P(X_{H} = y_H | X_T = x_T)\nonumber
\\[3pt]
&&\hspace*{110pt}{} \times\bigl\{P(X_w = 1, X_{{H^*}\setminus\{w\}} = y_{{H^*}\setminus
\{w\}} |X_{T^*} = x_{T^*})\nonumber
\\[3pt]
&&\hspace*{127pt}{} + P(X_w = 0, X_{{H^*}\setminus\{w\}} = y_{{H^*}\setminus\{w\}} |
X_{T^*} = x_{T^*}) \bigr\}.\nonumber
\end{eqnarray}
Note that by the definition of reducibility, $D^* \setminus H^*
\subseteq O \cup(W \setminus B)$, so $D^* \setminus H^*$ does~not
contain any ``flipped'' vertices; hence, $D^* \cap B \subseteq H^*$.
Further, $D^* \subseteq O \cup W$.

By\vspace*{1pt} Proposition \ref{propanc-dist-part}, applied to $H^*$, $D^*$ and $O \cup W$,
$ [O \cup W]_{\G} = [(O \cup W)\setminus D^*]_{\G} \cup[D^*]_{\G}$.
Thus, every head
$H^\dagger\in[O \cup W]_{\G}$ which contains a vertex in \mbox{$D^*
\setminus H^*$} is
such that $H^\dagger\subseteq D^*$. Hence, by applying
Lemma \ref{lemords}
to $D^*$,
it follows that $H^\dagger\prec H^*$ [since
$H^\dagger\subseteq D^* \setminus H^* \subseteq\an_{\mathcal G}(H^*)$
rules out
$H^* \prec H^\dagger$, while $H^*,H^\dagger\subseteq D^*$ rules out (iii)].
Thus, $D^*$ is made up of $H^*$ and
the heads which precede it under~$\prec$, and hence also under~$\prec^*$.
%

Suppose we replace $[O \cup W]_\G$ with $[(O \cup W) \setminus S]_\G$
for some $S \subseteq
H^*$; from Lemma \ref{lemheadfirst}, it is clear that only heads
which precede $H^*$ under $\prec^*$ will be affected, so in particular:%
%
%
%
%
%
\begin{eqnarray} \label{eqnpartition}
\bigl[(O\cup W) \setminus H^*\bigr]_{\G} &=& \bigl[(O\cup W)\setminus
D^*\bigr]_{\G} \cup\bigl[D^* \setminus H^*\bigr]_{\G}\quad\mbox{and}
\nonumber\\[-7pt]\\[-7pt]
\bigl[(O\cup W)\setminus\{w\}\bigr]_{\G} &=& \bigl[(O
\cup W)\setminus D^*\bigr]_{\G} \cup\bigl[D^*\setminus\{w\}
\bigr]_{\G}.\nonumber
\end{eqnarray}
It follows that to establish (\ref{eqnsuffice1}) it suffices to show:
%
\begin{eqnarray}\label{eqnsuffice2}
&& \prod_{H \in[D^*\setminus\{w\}]_{\G}} P(X_{H} =
y_H | X_T = x_T) \nonumber
\\[3pt]
&&\qquad = \bigl\{P(X_w = 0, X_{H^*\setminus\{w\}} = y_{H^* \setminus\{w\}} |
X_{T^*} = x_{T^*})
\nonumber\\[-7pt]\\[-7pt]
&&\hspace*{3pt}\quad\qquad{}+ P(X_w = 1, X_{H^*\setminus\{w\}} = y_{H^* \setminus\{w\}} |
X_{T^*} = x_{T^*}) \bigr\}
\nonumber
\\[3pt]
&&\quad\qquad{}\times\prod_{H\in[D^*\setminus H^*]_{\G}} P(X_{H} =
y_H | X_T = x_T).\nonumber
\end{eqnarray}
Let $z_{D^* \setminus\{w\}} \equiv y_{D^*\setminus\{w\}}$ and
$z_{V\setminus D^*} \equiv x_{V\setminus D^*}$ (with $z_w$ remaining free).
Since $D^* \cap B \subseteq H^*$, applying Lemma \ref{lemrmv-one} to
$D^*$ and $w$ using the values of $z_V$ gives
\begin{eqnarray*}
&& P(X_w = z_w | X_{(H^* \cup T^*) \setminus\{w\}} =
z_{(H^* \cup T^*) \setminus\{w\}}) \prod_{H \in[D^* \setminus\{w\}
]_\G} P(X_H =
z_H | X_T = z_T)
\\
&&\qquad = \prod_{H \in[D^*]_\G} P(X_H =
z_H | X_T = z_T)
\\
&&\qquad = P(X_{H^*} = z_{H^*} | X_{T^*} =
z_{T^*}) \prod_{H \in[D^* \setminus H^*]_\G} P(X_H =
z_H | X_T = z_T).
\end{eqnarray*}
Summing both sides of the equation over $z_w$ yields
(\ref{eqnsuffice2}). Thus, (\ref{eqngsum}) holds.

It remains to demonstrate that if $B \cup\{w\} \neq W$, the triples
$(x_V, B \cup\{w\}, W)$ and $(x_V, B, W \setminus\{w\})$ are also
reducible.

For the
first, consider $H \in[O \cup W]_\G$ with $H \cap(W \setminus(B \cup
\{w\})) \neq\varnothing$. Let $D \equiv\dis_{\an(H)}(H) \subseteq
O\cup W$;
by construction $D$ has ancestrally closed districts.
Since $H \cap(W \setminus B) \supseteq H \cap(W \setminus(B \cup
\{w\})) \neq\varnothing$, by the reducibility of $(x_V, B, W)$,
$D \setminus H \subseteq O \cup(W \setminus B)$.
It is sufficient to show that $w \notin D \setminus H$.
Since by Proposition~\ref{propanc-dist-part},
$ [O \cup W]_{\G} = [(O \cup W)\setminus D]_{\G} \cup[D]_{\G}$,
if $w \in D \cap H^*$ then $H^* \in[D]_{\mathcal G}$. If $H = H^*$, then
$w \notin D \setminus H$. If $H \neq H^*$, then applying
Lemma \ref{lemords} we have $H^* \prec H$ (by the same argument as above).
But this contradicts that $H^*$ is a maximal head in $[O\cup W]_{\G}$
such that
$H^* \cap(W \setminus B) \neq\varnothing$.
Hence, $(x_V, B \cup\{w\}, W)$ is reducible.
%

We now consider $(x_V, B, W \setminus\{w\})$. Let $H \in[O
\cup(W \setminus\{w\})]$, with $H \cap((W \setminus
\{w\}) \setminus B) \neq\varnothing$. Again, let $D \equiv\dis_{\an(H)}(H)$.

First suppose $H \in[(O \cup W) \setminus D^*]_\G$ then, by (\ref
{eqnpartition}), $H \in[O \cup W]_\G$.
We showed above that if $H \in[O \cup W]_\G$ and $H \cap(W
\setminus(B \cup\{w\}))\neq\varnothing$ then
$D \setminus H \subseteq(W \setminus(B \cup\{w\}))$. This is
sufficient since $W \setminus(B \cup\{w\}) =
(W\setminus\{w\}) \setminus B$.

If $H \notin[(O \cup W) \setminus D^*]_\G$ then (\ref{eqnpartition}) implies
$H \in[D^* \setminus\{w\}]_\G$. Lemma \ref{lemancdist} applied to
$D^*$ implies that
$D^*\setminus\{w\}$ has ancestrally closed districts, so $D \subseteq
D^*\setminus\{w\}$.
Since $D \subseteq D^*$, if a vertex $v$ is not barren in $D$ then $v
\notin\barren_\G(D^*) = H^*$.
Hence, $H^* \cap D \subseteq \barren_\G(D) = H$.
%
Thus,
\[
D \setminus H \subseteq D \setminus H^* \subseteq D^* \setminus H^*
\subseteq O
\cup(W \setminus B),
\]
where the third inclusion follows from the reducibility of $(x_V, B,
W)$ and the choice of $H^*$. But since $D \subseteq D^*\setminus\{w\}
$, we have
$D \setminus H \subseteq
O \cup((W \setminus\{w\}) \setminus B)$ as required.
\end{pf*}
\end{appendix}


%

\printaddresses

\begin{thebibliography}{29}
\bibitem[\protect\citeauthoryear{Andersson, Madigan and Perlman}{2001}]{amp}
%
\begin{barticle}[mr]
\bauthor{\bsnm{Andersson},~\bfnm{Steen~A.}\binits{S.~A.}},
\bauthor{\bsnm{Madigan},~\bfnm{David}\binits{D.}} \AND
\bauthor{\bsnm{Perlman},~\bfnm{Michael~D.}\binits{M.~D.}}
(\byear{2001}).
\btitle{Alternative {M}arkov properties for chain graphs}.
\bjournal{Scand. J. Stat.}
\bvolume{28}
\bpages{33--85}.
\bid{doi={10.1111/1467-9469.00224}, issn={0303-6898}, mr={1844349}}
\end{barticle}
%
\bptok{imsref}%
\endbibitem

\bibitem[\protect\citeauthoryear{Bartolucci, Colombi and Forcina}{2007}]{bart07}
%
\begin{barticle}[mr]
\bauthor{\bsnm{Bartolucci},~\bfnm{Francesco}\binits{F.}},
\bauthor{\bsnm{Colombi},~\bfnm{Roberto}\binits{R.}} \AND
\bauthor{\bsnm{Forcina},~\bfnm{Antonio}\binits{A.}}
(\byear{2007}).
\btitle{An extended class of marginal link functions for modelling
contingency tables by equality and inequality constraints}.
\bjournal{Statist. Sinica}
\bvolume{17}
\bpages{691--711}.
\bid{issn={1017-0405}, mr={2398430}}
\end{barticle}
%
\bptok{imsref}%
\endbibitem

\bibitem[\protect\citeauthoryear{Bergsma and Rudas}{2002}]{br02}
%
\begin{barticle}[mr]
\bauthor{\bsnm{Bergsma},~\bfnm{Wicher~P.}\binits{W.~P.}} \AND
\bauthor{\bsnm{Rudas},~\bfnm{Tam{\'a}s}\binits{T.}}
(\byear{2002}).
\btitle{Marginal models for categorical data}.
\bjournal{Ann. Statist.}
\bvolume{30}
\bpages{140--159}.
\bid{doi={10.1214/aos/1015362188}, issn={0090-5364}, mr={1892659}}
\end{barticle}
%
\bptok{imsref}%
\endbibitem

\bibitem[\protect\citeauthoryear{Dawid}{1979}]{dawidcondind}
%
\begin{barticle}[mr]
\bauthor{\bsnm{Dawid},~\bfnm{A.~P.}\binits{A.~P.}}
(\byear{1979}).
\btitle{Conditional independence in statistical theory}.
\bjournal{J. Roy. Statist. Soc. Ser. B}
\bvolume{41}
\bpages{1--31}.
\bid{issn={0035-9246}, mr={0535541}}
\bptnote{check related}%
\end{barticle}
%
\bptok{imsref}%
\endbibitem

\bibitem[\protect\citeauthoryear{Dawid and Didelez}{2010}]{dawiddidelez2010}
%
\begin{barticle}[mr]
\bauthor{\bsnm{Dawid},~\bfnm{A.~Philip}\binits{A.~P.}} \AND
\bauthor{\bsnm{Didelez},~\bfnm{Vanessa}\binits{V.}}
(\byear{2010}).
\btitle{Identifying the consequences of dynamic treatment strategies: A
decision-theoretic overview}.
\bjournal{Stat. Surv.}
\bvolume{4}
\bpages{184--231}.
\bid{doi={10.1214/10-SS081}, issn={1935-7516}, mr={2740837}}
\end{barticle}
%
\bptok{imsref}%
\endbibitem

\bibitem[\protect\citeauthoryear{Drton}{2009}]{drton09}
%
\begin{barticle}[mr]
\bauthor{\bsnm{Drton},~\bfnm{Mathias}\binits{M.}}
(\byear{2009}).
\btitle{Discrete chain graph models}.
\bjournal{Bernoulli}
\bvolume{15}
\bpages{736--753}.
\bid{doi={10.3150/08-BEJ172}, issn={1350-7265}, mr={2555197}}
\end{barticle}
%
\bptok{imsref}%
\endbibitem

\bibitem[\protect\citeauthoryear{Drton and Xiao}{2010}]{drton10}
%
\begin{bincollection}[mr]
\bauthor{\bsnm{Drton},~\bfnm{Mathias}\binits{M.}} \AND
\bauthor{\bsnm{Xiao},~\bfnm{Han}\binits{H.}}
(\byear{2010}).
\btitle{Smoothness of {G}aussian conditional independence models}.
In \bbooktitle{Algebraic Methods in Statistics and Probability {II}}.
\bseries{Contemp. Math.}
\bvolume{516}
\bpages{155--177}.
\bpublisher{Amer. Math. Soc.},
\blocation{Providence, RI}.
\bid{doi={10.1090/conm/516/10173}, mr={2730747}}
\end{bincollection}
%
\bptok{imsref}%
\endbibitem

\bibitem[\protect\citeauthoryear{Evans and Richardson}{2010}]{evans10}
%
\begin{bincollection}[auto:STB|2014/02/12|14:17:21]
\bauthor{\bsnm{Evans},~\bfnm{R.~J.}\binits{R.~J.}} \AND
\bauthor{\bsnm{Richardson},~\bfnm{T.~S.}\binits{T.~S.}}
(\byear{2010}).
\btitle{Maximum likelihood fitting of acyclic directed mixed graphs to
binary data}.
In \bbooktitle{Proceedings of the 26th Conference on Uncertainty in
Artificial Intelligence}
\bpages{177--184}.
\bpublisher{AUAI Press}, \blocation{Corvallis, OR}.
\end{bincollection}
%
\bptok{imsref}%
\endbibitem

\bibitem[\protect\citeauthoryear{Evans and
Richardson}{2013}]{evansrichardson2013}
%
\begin{barticle}[mr]
\bauthor{\bsnm{Evans},~\bfnm{Robin~J.}\binits{R.~J.}} \AND
\bauthor{\bsnm{Richardson},~\bfnm{Thomas~S.}\binits{T.~S.}}
(\byear{2013}).
\btitle{Marginal log-linear parameters for graphical {M}arkov models}.
\bjournal{J. R. Stat. Soc. Ser. B Stat. Methodol.}
\bvolume{75}
\bpages{743--768}.
\bid{doi={10.1111/rssb.12020}, issn={1369-7412}, mr={3091657}}
\end{barticle}
%
\bptok{imsref}%
\endbibitem

\bibitem[\protect\citeauthoryear{Hammersley and Clifford}{1971}]{hammersley71}
%
\begin{bmisc}[auto:STB|2014/02/12|14:17:21]
\bauthor{\bsnm{Hammersley},~\bfnm{J.~M.}\binits{J.~M.}} \AND
\bauthor{\bsnm{Clifford},~\bfnm{P.}\binits{P.}}
(\byear{1971}).
\bhowpublished{Markov fields on finite graphs and lattices. Unpublished
manuscript}.
\end{bmisc}
%
\bptok{imsref}%
\endbibitem

\bibitem[\protect\citeauthoryear{Huang and Valtorta}{2006}]{huang06do}
%
\begin{bincollection}[auto:STB|2014/02/12|14:17:21]
\bauthor{\bsnm{Huang},~\bfnm{Y.}\binits{Y.}} \AND
\bauthor{\bsnm{Valtorta},~\bfnm{M.}\binits{M.}}
(\byear{2006}).
\btitle{Pearl's calculus of interventions is complete}.
In \bbooktitle{Proceedings of the 22nd Conference On Uncertainty in
Artificial Intelligence}.
\bpublisher{AUAI Press}, \blocation{Arlington, VA}.
\end{bincollection}
%
\bptok{imsref}%
\endbibitem

\bibitem[\protect\citeauthoryear{Kass and Vos}{1997}]{kass97}
%
\begin{bbook}[mr]
\bauthor{\bsnm{Kass},~\bfnm{Robert~E.}\binits{R.~E.}} \AND
\bauthor{\bsnm{Vos},~\bfnm{Paul~W.}\binits{P.~W.}}
(\byear{1997}).
\btitle{Geometrical Foundations of Asymptotic Inference}.
\bpublisher{Wiley},
\blocation{New York}.
\bid{doi={10.1002/9781118165980}, mr={1461540}}
\end{bbook}
%
\bptok{imsref}%
\endbibitem

\bibitem[\protect\citeauthoryear{Lauritzen}{1996}]{lau96}
%
\begin{bbook}[mr]
\bauthor{\bsnm{Lauritzen},~\bfnm{Steffen~L.}\binits{S.~L.}}
(\byear{1996}).
\btitle{Graphical Models}.
\bpublisher{Oxford Univ. Press},
\blocation{New York}.
\bid{mr={1419991}}
\end{bbook}
%
\bptok{imsref}%
\endbibitem

\bibitem[\protect\citeauthoryear{Pearl}{1988}]{pearl88}
%
\begin{bbook}[mr]
\bauthor{\bsnm{Pearl},~\bfnm{Judea}\binits{J.}}
(\byear{1988}).
\btitle{Probabilistic Reasoning in Intelligent Systems: Networks of
Plausible Inference}.
\bpublisher{Morgan Kaufmann},
\blocation{San Mateo, CA}.
\bid{mr={0965765}}
\end{bbook}
%
\bptok{imsref}%
\endbibitem

\bibitem[\protect\citeauthoryear{Pearl}{1995}]{pearlbiom}
%
\begin{barticle}[mr]
\bauthor{\bsnm{Pearl},~\bfnm{Judea}\binits{J.}}
(\byear{1995}).
\btitle{Causal diagrams for empirical research}.
\bjournal{Biometrika}
\bvolume{82}
\bpages{669--710}.
\bid{doi={10.1093/biomet/82.4.669}, issn={0006-3444}, mr={1380809}}
\end{barticle}
%
\bptok{imsref}%
\endbibitem

\bibitem[\protect\citeauthoryear{Pearl}{2009}]{pearl2009}
%
\begin{bbook}[mr]
\bauthor{\bsnm{Pearl},~\bfnm{Judea}\binits{J.}}
(\byear{2009}).
\btitle{Causality},
\bedition{2nd} ed.
\bpublisher{Cambridge Univ. Press},
\blocation{Cambridge}.
\bid{doi={10.1017/CBO9780511803161}, mr={2548166}}
\end{bbook}
%
\bptok{imsref}%
\endbibitem

\bibitem[\protect\citeauthoryear{Pearl and Robins}{1995}]{pearlrobins1995}
%
\begin{bincollection}[mr]
\bauthor{\bsnm{Pearl},~\bfnm{Judea}\binits{J.}} \AND
\bauthor{\bsnm{Robins},~\bfnm{James}\binits{J.}}
(\byear{1995}).
\btitle{Probabilistic evaluation of sequential plans from causal models
with hidden variables}.
In \bbooktitle{Uncertainty in Artificial Intelligence ({M}ontreal,
{PQ}, 1995)}
\bpages{444--453}.
\bpublisher{Morgan Kaufmann},
\blocation{San Francisco, CA}.
\bid{mr={1615028}}
\end{bincollection}
%
\bptok{imsref}%
\endbibitem

\bibitem[\protect\citeauthoryear{Richardson}{2003}]{richardson03}
%
\begin{barticle}[mr]
\bauthor{\bsnm{Richardson},~\bfnm{Thomas}\binits{T.}}
(\byear{2003}).
\btitle{Markov properties for acyclic directed mixed graphs}.
\bjournal{Scand. J. Stat.}
\bvolume{30}
\bpages{145--157}.
\bid{doi={10.1111/1467-9469.00323}, issn={0303-6898}, mr={1963898}}
\end{barticle}
%
\bptok{imsref}%
\endbibitem

\bibitem[\protect\citeauthoryear{Richardson}{2009}]{richardson09}
%
\begin{bincollection}[auto:STB|2014/02/12|14:17:21]
\bauthor{\bsnm{Richardson},~\bfnm{T.~S.}\binits{T.~S.}}
(\byear{2009}).
\btitle{A factorization criterion for acyclic directed mixed graphs}.
In \bbooktitle{Proceedings of the 25th Conference on Uncertainty in
Artificial Intelligence}
\bpages{462--470}.
\bpublisher{AUAI Press}, \blocation{Arlington, VA}.
\end{bincollection}
%
\bptok{imsref}%
\endbibitem

\bibitem[\protect\citeauthoryear{Richardson and Spirtes}{2002}]{richardson02}
%
\begin{barticle}[mr]
\bauthor{\bsnm{Richardson},~\bfnm{Thomas}\binits{T.}} \AND
\bauthor{\bsnm{Spirtes},~\bfnm{Peter}\binits{P.}}
(\byear{2002}).
\btitle{Ancestral graph {M}arkov models}.
\bjournal{Ann. Statist.}
\bvolume{30}
\bpages{962--1030}.
\bid{doi={10.1214/aos/1031689015}, issn={0090-5364}, mr={1926166}}
\end{barticle}
%
\bptok{imsref}%
\endbibitem

\bibitem[\protect\citeauthoryear{Robins and Richardson}{2011}]{robinsmcm2011}
%
\begin{bincollection}[auto:STB|2014/02/12|14:17:21]
\bauthor{\bsnm{Robins},~\bfnm{J.~M.}\binits{J.~M.}} \AND
\bauthor{\bsnm{Richardson},~\bfnm{T.~S.}\binits{T.~S.}}
(\byear{2011}).
\btitle{Alternative graphical causal models and the identification of
direct effects}.
In \bbooktitle{Causality and Psychopathology: Finding the Determinants
of Disorders and Their Cures}
(\beditor{\bfnm{Patrick}\binits{P.}~\bsnm{Shrout}},
\beditor{\bfnm{Katherine}\binits{K.}~\bsnm{Keyes}} \AND
\beditor{\bfnm{Katherine}\binits{K.}~\bsnm{Ornstein}}, eds.)
\bvolume{6}
\bpages{1--52}.
\bpublisher{Oxford Univ. Press},
\blocation{London}.
\end{bincollection}
%
\bptok{imsref}%
\endbibitem

\bibitem[\protect\citeauthoryear{Sadeghi}{2013}]{sadeghi11b}
%
\begin{barticle}[mr]
\bauthor{\bsnm{Sadeghi},~\bfnm{Kayvan}\binits{K.}}
(\byear{2013}).
\btitle{Stable mixed graphs}.
\bjournal{Bernoulli}
\bvolume{19}
\bpages{2330--2358}.
\bid{doi={10.3150/12-BEJ454}, issn={1350-7265}, mr={3160556}}
\end{barticle}
%
\bptok{imsref}%
\endbibitem

\bibitem[\protect\citeauthoryear{Sadeghi and Lauritzen}{2014}]{sadeghi11}
%
\begin{barticle}[mr]
\bauthor{\bsnm{Sadeghi},~\bfnm{Kayvan}\binits{K.}} \AND
\bauthor{\bsnm{Lauritzen},~\bfnm{Steffen}\binits{S.}}
(\byear{2014}).
\btitle{Markov properties for mixed graphs}.
\bjournal{Bernoulli}
\bvolume{20}
\bpages{676--696}.
\bid{doi={10.3150/12-BEJ502}, issn={1350-7265}, mr={3178514}}
\end{barticle}
%
\bptok{imsref}%
\endbibitem

\bibitem[\protect\citeauthoryear{Shpitser and Pearl}{2006a}]{shpitser06id}
%
\begin{bincollection}[auto:STB|2014/02/12|14:17:21]
\bauthor{\bsnm{Shpitser},~\bfnm{I.}\binits{I.}} \AND
\bauthor{\bsnm{Pearl},~\bfnm{J.}\binits{J.}}
(\byear{2006}a).
\btitle{Identification of joint interventional distributions in
recursive semi-Markovian causal models}.
In \bbooktitle{Proceedings of the 21st National Conference on
Artificial Intelligence}.
\bpublisher{AAAI Press}, \blocation{Menlo Park, CA}.
\end{bincollection}
%
\bptok{imsref}%
\endbibitem

\bibitem[\protect\citeauthoryear{Shpitser and Pearl}{2006b}]{shpitser06idc}
%
\begin{bincollection}[auto:STB|2014/02/12|14:17:21]
\bauthor{\bsnm{Shpitser},~\bfnm{I.}\binits{I.}} \AND
\bauthor{\bsnm{Pearl},~\bfnm{J.}\binits{J.}}
(\byear{2006}b).
\btitle{Identification of conditional interventional distributions}.
In \bbooktitle{Proceedings of the 22nd Conference on Uncertainty in
Artificial Intelligence}
\bpages{437--444}.
\bpublisher{AUAI Press}, \blocation{Arlington, VA}.
\end{bincollection}
%
\bptok{imsref}%
\endbibitem

\bibitem[\protect\citeauthoryear{Silva and Ghahramani}{2009}]{silva09}
%
\begin{barticle}[mr]
\bauthor{\bsnm{Silva},~\bfnm{Ricardo}\binits{R.}} \AND
\bauthor{\bsnm{Ghahramani},~\bfnm{Zoubin}\binits{Z.}}
(\byear{2009}).
\btitle{The hidden life of latent variables: {B}ayesian learning with
mixed graph models}.
\bjournal{J. Mach. Learn. Res.}
\bvolume{10}
\bpages{1187--1238}.
\bid{issn={1532-4435}, mr={2520804}}
\end{barticle}
%
\bptok{imsref}%
\endbibitem

\bibitem[\protect\citeauthoryear{Spirtes, Glymour and Scheines}{1993}]{cps93}
%
\begin{bbook}[mr]
\bauthor{\bsnm{Spirtes},~\bfnm{Peter}\binits{P.}},
\bauthor{\bsnm{Glymour},~\bfnm{Clark}\binits{C.}} \AND
\bauthor{\bsnm{Scheines},~\bfnm{Richard}\binits{R.}}
(\byear{1993}).
\btitle{Causation, Prediction, and Search}.
\bseries{Lecture Notes in Statistics}
\bvolume{81}.
\bpublisher{Springer},
\blocation{New York}.
\bid{doi={10.1007/978-1-4612-2748-9}, mr={1227558}}
\end{bbook}
%
\bptok{imsref}%
\endbibitem

\bibitem[\protect\citeauthoryear{Tian and Pearl}{2002}]{tian02b}
%
\begin{bincollection}[auto:STB|2014/02/12|14:17:21]
\bauthor{\bsnm{Tian},~\bfnm{J.}\binits{J.}} \AND
\bauthor{\bsnm{Pearl},~\bfnm{J.}\binits{J.}}
(\byear{2002}).
\btitle{A general identification condition for causal effects}.
In \bbooktitle{Proceedings of the 18th National Conference on
Artificial Intelligence}.
\bpublisher{AAAI Press},
\blocation{Menlo Park}.
\end{bincollection}
%
\bptok{imsref}%
\endbibitem

\bibitem[\protect\citeauthoryear{Wermuth}{2011}]{wermuth11}
%
\begin{barticle}[mr]
\bauthor{\bsnm{Wermuth},~\bfnm{Nanny}\binits{N.}}
(\byear{2011}).
\btitle{Probability distributions with summary graph structure}.
\bjournal{Bernoulli}
\bvolume{17}
\bpages{845--879}.
\bid{doi={10.3150/10-BEJ309}, issn={1350-7265}, mr={2817608}}
\end{barticle}
%
\bptok{imsref}%
\endbibitem

\end{thebibliography}
\end{document}